\documentclass[a4paper,oneside,11pt]{amsart}

\usepackage{enumerate,stackrel}

\usepackage[english]{babel}
\usepackage{fancyhdr}
\usepackage{amsmath}
\usepackage{mathrsfs}
\usepackage{color}

\usepackage{comment}
\let\oldmarginpar\marginpar
\renewcommand\marginpar[1]{\-\oldmarginpar[\raggedleft\small\sf
#1] {\raggedright\small\sf #1}}

\usepackage{amssymb}

\newtheorem{definition}{Definition}[section]

\newtheorem{theorem}[definition]{Theorem}

\newtheorem{corollary}[definition]{Corollary}
\newtheorem{lemma}[definition]{Lemma}
\newtheorem{remark}[definition]{Remark}

\numberwithin{equation}{section}

\def\bkR{{\rm I\kern-.17em R}}

\def\bkN{{\rm I\kern-.17em N}}

\def\bkC{{\rm \kern.24em \vrule width.05em height1.4ex depth-.05ex
     \kern-.26em C}}

\def\bkK{{\rm I\kern-.22em K}}
\def\bkP{{\rm I\kern-.22em P}}

\def\ds{\displaystyle}

\def\cB{{\mcal B}}

\newcommand{\mcal}{\mathcal}
\newcommand{\balpha}{\mbox{\boldmath $\alpha$\unboldmath}}

\newcommand{\bell}{\mbox{\boldmath $\ell$\unboldmath}}

\newcommand{\bpr}{\begin{proof}}
\newcommand{\epr}{\hfill $\square$ \end{proof}}
\newcommand{\bprof}{\begin{proofof}}
\newcommand{\eprof}{\hfill $\square$ \end{proofof}}
\newcommand{\bproff}{\begin{proofoff}}
\newcommand{\eproff}{\hfill $\square$ \end{proofoff}}
\begin{document}


\title{Optimal local embeddings of Besov spaces involving only slowly varying smoothness}

\author{J\'ulio S. Neves}
\author{Bohum\'{\i}r Opic}

\address{J\'ulio Severino Neves, CMUC, Department of Mathematics, University of Coimbra, Apartado 3008, 3001-454 Coimbra, Portugal}
\email{jsn@mat.uc.pt}

\address{Bohum\'{\i}r Opic, Department of Mathematical Analysis, 
Faculty of Mathematics and Physics, Charles University,
 Sokolovsk\'a 83, 186 75 Prague 8, Czech Republic}
\email{opic@karlin.mff.cuni.cz}

\keywords{Besov spaces involving only a slowly varying smoothness; sharp and optimal
embeddings, limiting real interpolation, weighted inequalities}

\subjclass[2000]{46E35, 46E30, 26A15, 26A12, 46B70, 26D10, 26D15}

\begin{abstract}
The aim of the paper is to establish (local) optimal embeddings of Besov spaces $B^{0,b}_{p,r}$ involving only a slowly
varying smoothness $b$. In general, our target spaces are outside of the scale of
Lorentz-Karamata spaces and are related to small Lebesgue spaces. In particular, we improve
results from \cite{CGO2}, where the targets are (local) Lorentz-Karamata
spaces. To derive such results, we apply limiting real interpolation techniques
and weighted Hardy-type inequalities.
\end{abstract}

\thanks{The research has been supported by the grant P201-18-00580S of the Grant
Agency of the Czech Republic and by Centre of Mathematics of the University of Coimbra -- UID/MAT/00324/2019, funded by the Portuguese Government through FCT/MEC and co-funded by the European Regional Development Fund through the Partnership Agreement PT2020. } 

\maketitle


\section{Introduction}

Classical Besov spaces $B^{s,b}_{p,r}$ play an important role in numerous parts of mathematics. However, it has gradually become clear that, to solve some problems, 
Besov spaces with smoothness which can be more finely turned, i.e., Besov spaces with generalized smoothness, are essential. 
These spaces have been studied especially by the Soviet mathematical school (cf. \cite [Section 8]{KL87}) and a~lot of attention has been paid to optimal embeddings and to growth envelopes of such spaces (see, e.g., \cite{Gol-Russian}, \cite{GoKe}, \cite{Net89}, 
\cite{FaLe06}, 
\cite{CaFa04}, 
\cite{CM04a},\cite{CM04b}, \cite{CaHa04}, 
\cite{GuOp04},  
\cite{Moura}, 
\cite{GO07:SEBtS}, 
\cite[Chapter~1]{Tri06}, \cite{Gol07}, \cite{CGO1}, \cite{CGO2}, 
etc.)

Besov spaces  $B^{0,b}_{p,r}$ involving the zero classical smoothness and  logarithmic smoothness $b$ are particular cases of Besov spaces 
with generalized smoothness. 
They appear (probably for the first time) already in  \cite{DeVRS} in a connection with the weak-type interpolation. During the recent years these spaces have attracted an increasing attention (see \cite{CGO1}, 
\cite{CD1}, \cite{GOTT2014}, \cite{CDSFM2015}, \cite{CD2}, 
\cite{CD3}, 
\cite{CDT}, 
\cite{D2016},  \cite{D2017}, \cite{BeCo}, \cite{CDK2018}, etc.). The problem of optimal embeddings in this setting is much harder to study.


In \cite{CGO2} the authors have characterized 
(with easily verifiable conditions) local embeddings of Besov spaces $B^{0,b}_{p,r}$  
involving the zero classical smoothness and involving only a slowly varying smoothness $b$ into classical Lorentz spaces.  These results have been then applied to establish sharp local embeddings of Besov spaces in question into Lorentz-Karamata spaces and to determinate the growth envepoles of spaces $B^{0,b}_{p,r}$. In particular, the following three theorems (adapted to our notation, cf. Remarks \ref{remark:considerationsDefinBesov} and \ref{remark:13} below) have been proved there.\footnotemark$^)$
\footnotetext{$^{)}$
Besov spaces are defined by means of the modulus of continuity, we refer to Section 2 for precise definition.}

  \begin{theorem}[{\cite[Theorem 3.2 (i)]{CGO2}}]\label{T1:CGO2}
Let $n\in\bkN$, $1 \leq p <+ \infty$, $1 \leq r \le  q \leq +\infty$. Let $b$ be a slowly varying function on the interval $(0, +\infty)$ \rm{(}notation $b\in SV(0,+\infty)$\rm{)} satisfying 
\begin{equation}
	\label{(1*)}
\|t^{-1/r}\;b(t)\|_{r;(0,1)}=+\infty\quad\text{and}\quad
\|t^{-1/r}\;b(t)\|_{r;(1,+\infty)}<+\infty,
\end{equation} 
and let $b_{r,n}$ be defined by 
\begin{equation}
\label{b_r}
b_{r,n}(t) := \| s^{-1/r} b(s^{1/n}) \|_{r;(t,+\infty)}\quad\text{for all }\quad t\in(0,+\infty).
\end{equation}
Put $\rho=+\infty$ if $p \le q$ and define $\rho$ by $\frac{1}{\rho}=\frac{1}{q}-\frac{1}{p}$ if $q<p$. 
Assume that $w$ is a non-negative measurable function on $(0,1)$ and
$${\mcal W}_q(t):=\| w(s)\|_{q;(0,t)}, \quad  t \in (0,1].$$
Then
\begin{equation}\label{T1'}
    \| w(t) f^*(t) \|_{q;(0,1)} \precsim \| f
    \|_{B^{0,b}_{p,r}}\quad\text{for all}\quad f \in B^{0,b}_{p,r}(\bkR^n)
\end{equation}
if and only if  
\begin{equation}\label{T1.1}
{\mcal W}_q(1)+\|s^{-\frac{1}{p}-\frac{1}{\rho}}{\mcal W}_q(s)\|_{\rho;(t,1)}\precsim b_{r,n}(t)\quad \text{for all}\quad t\in (0,1).
\end{equation}
\end{theorem}

\medskip

\begin{theorem}[{\cite[Theorem 3.3]{CGO2}}]
\label{main*:CGO2}
 Let $n\in\bkN$, $1 \leq p < +\infty$, $\, 1 \leq r \leq+ \infty$, $0 < q \leq +\infty$ 
and let $b\in SV(0,+\infty)$ satisfy \eqref{(1*)}.
Let $b_{r,n}$ be giving by
\eqref{b_r} and 
define, for all $t \in (0,+\infty)$, 
\begin{equation}
\label{O61*}
	\tilde{b}(t):= \left\{
\begin{array}{ll}
(b_{r,n}(t))^{1-r/q+r/\max\{ p,q\}} (b(t^{1/n}))^{r/q-r/\max\{ p,q\}} & \mbox{ if }\, r \not= +\infty \\
b_{\infty,n}(t) & \mbox{ if }\, r = +\infty
\end{array}
\right..
\end{equation}
Then 
\begin{equation}
\label{embedding*}
	\| t^{1/p-1/q} \tilde{b}(t) f^*(t) \|_{q;(0,1)} \precsim \| f \|_{B^{0,b}_{p,r}}\quad\text{for all}\quad f \in B^{0,b}_{p,r}(\bkR^n)
\end{equation}
%
if and only if $q \geq r$.
\end{theorem}
\medskip

\begin{theorem}[{\cite[Theorem 3.4]{CGO2}}]
\label{kappa*:CGO2}
Let $n\in\bkN$, $1 \leq p <+ \infty$, $1 \leq r \leq q \leq +\infty$ and let 
$\, b \in SV(0,+\infty)\,$ satisfy~\eqref{(1*)}.
Define $b_{r,n}$ and $\tilde{b}$  
by \eqref{b_r} and \eqref{O61*}, respectively. Let $\kappa$ be a~non-negative and non-increasing function on $(0, 1)$. Then
\begin{equation}
\label{kappa_embedding*}
	\| t^{1/p-1/q} \tilde{b}(t) \kappa(t) f^*(t) \|_{q;(0,1)} \precsim \| f \|_{B^{0,b}_{p,r}}\quad\text{for all}\quad f \in B^{0,b}_{p,r}(\bkR^n)
\end{equation}
if and only if $\kappa$ is bounded.
\end{theorem}

Theorems  \ref{main*:CGO2} and \ref{kappa*:CGO2} describe the sharp continuous embeddings of the Besov space  $B^{0,b}_{p,r}(\bkR^n)$ into the local Lorentz-Karamata space $L_{p, q; \tilde{b}}^{loc}(\bkR^n)$ (we refer to Section \ref{section2} for the precise definition of spaces in question).
Namely, these theorems imply that
 \begin{equation} \label{vnorLocal}
B^{0,b}_{p,r}(\bkR^n)\hookrightarrow L_{p, q; \tilde{b}}^{loc}(\bkR^n)
\end{equation} 
and that this embedding is sharp within the scale of local Lorentz-Karamata spaces.

For  characterization of compact embeddings of spaces $B^{0,b}_{p,r}(\bkR^n)$ into Lorentz-Karamata spaces we refer to \cite{CGO3}.

The aim of this paper is to improve embedding \eqref{vnorLocal}. Namely, we are going to prove (cf. Theorems \ref{theo:mainresult} and \ref{theorem24} of Section \ref{sectionmain}) that the target space in \eqref{vnorLocal} can be replaced by the better space
$Z(\bkR^n)=Z_{p,q,n,\bar{b}}(\bkR^n)$  if $p\neq q$, where $Z_{p,q,n,\bar{b}}(\bkR^n)$ is the set of all measurable $f$ functions on $\bkR^n$ satisfying 
\begin{equation}\label{defi:quasinormZspace}
\|f\|_{Z}=\|f\|_{Z_{p,q,n,\bar{b}}}:=\big\|t^{-1/q}\bar{b}(t^{1/n})\,\|f^*\|_{p;(0,t)}\big\|_{q;(0,+\infty)}<+\infty,
\end{equation}
with $\bar{b}=\bar{b}_{(r,q)}$ given by \eqref{deffunct:maintheo:1} below .

Moreover,  the target space in \eqref{vnorLocal} can still be replaced by $Z^{loc}(\bkR^n)=Z^{loc}_{p,q,n,\bar{b}}(\bkR^n)$, if $p\neq q$, which is a local version of $Z(\bkR^n)$, defined by replacing the interval $(0,+\infty)$ by $(0,1)$ in the (quasi-)norm \eqref{defi:quasinormZspace}, that is, $Z^{loc}_{p,q,n,\bar{b}}(\bkR^n)$ is the set of all measurable $f$ functions on $\bkR^n$ satisfying 
$$
\|f\|_{Z^{loc}}=\|f\|_{Z^{loc}_{p,q,n,\bar{b}}}:=\big\|t^{-1/q}\bar{b}(t^{1/n})\,\|f^*\|_{p;(0,t)}\big\|_{q;(0,1)}<+\infty,
$$
with $\bar{b}=\bar{b}_{(r,q)}$ given by \eqref{deffunct:maintheo:1} below .

Furthermore, we show (cf. Theorems \ref{theo:relatedtolemma14} and \ref{sharpnesswithrespecttofunction} below) that the embedding 
\begin{equation} \label{sharp}
B^{0,b}_{p,r}(\bkR^n)\hookrightarrow Z_{p,q,n,\bar{b}_{(r,q)}}(\bkR^n)
\end{equation}
is locally sharp with respect to the slowly varying function $\bar{b}_{(r,q)}$ and that embedding \eqref{sharp} with $q=r$ is optimal among all embeddings \eqref{sharp} 
with $1\le r\le q <+\infty$.

To prove our results, we make use of limiting interpolation, embedding theorems for spaces from two different scales of interpolation spaces (cf. Theorems \ref{Theorem5:circ},  
\ref{Theorem5:circ:STAR} below, which are of independent interest),  weighted inequalities, and some results from \cite{CGO2}.

In the particular case when $n=1$ and $b$ is of logarithmic type, the embedding of the subspace ${\mcal P}_{2\pi}(\bkR)\subset B_{p,r}^{0,b}(\bkR)$ of $2\pi$-periodical functions 
into the space $Z_{p,q,1,\bar{b}}(\bkR)$ follows from
\cite[Theorem 4.4]{COBOS201543}, where the authors used approximation spaces (including their limiting forms), reiteration of approximation constructions, the Nikolski\u{\i} inequality for trigonometric polynomials and limiting interpolation.

As in \cite{CGO2}, also in our paper  Besov spaces are defined by means of the modulus of continuity. Note that some authors use the Fourier-analytical approach to define Besov spaces with the zero classical smoothness and involving the logarithmic smoothnees $b$, where 
\begin{equation} \label{funkce b}
b(t):=(1+|\ln t|)^{\alpha}, \ \ t\in(0,+\infty), \quad\text{with a convenient \ }  \alpha \in \bkR. 
\end{equation}
However, if we denote the resulting space by $\cB^{0, b}_{p,r}(\bkR^n)$, with $b$ given by \eqref{funkce b}, then, in general, 
$$ B^{0,b}_{p,r}(\bkR^n)\neq \cB^{0, b}_{p,r}(\bkR^n)$$ (see \cite{SiTr},
\cite{CL2013}, \cite{CD2}).

When $b$ is given by \eqref{funkce b}, and $n\in\bkN$, then an analogue of embedding \eqref{sharp} with $B_{p,r}^{0,b}(\bkR^n)$ replaced by the subspace ${\mcal P}(\bkR^n)\subset B_{p,r}^{0,b}(\bkR^n)$ 
of periodical functions has been proved in \cite{D2016} (where also the corresponding result for the the subspace ${\mcal P}(\bkR^n)\subset \cB_{p,r}^{0,b}(\bkR^n)$ 
of periodical functions can be found, see also \cite{DomO17}).

The paper is organized as follows: Section 2 contains notation and preliminaries.
In Section~3 we present our main results on embeddings of Besov spaces involving only slowly varying smoothness, their sharpness and local optimality. In Section~4 we collect and prove some auxiliary results. The proofs of our main results from Section~3 are given in Sections~5, 6 and 7.


\section{Notation, definitions and basic properties}\label{section2}
\indent

As usual, $\bkR^n$ denotes the Euclidean $n$-dimensional space.
Throughout the paper $\mu_n$ is the $n$-dimensional Lebesgue measure
in $\bkR^n$ and $\Omega$ is a domain in $\bkR^n$.
We denote by $\chi_{\Omega}$ the characteristic function of $\Omega$
and put $|\Omega|_n=\mu_n(\Omega)$. The family of all extended
scalar-valued (real or complex) $\mu_n$-measurable functions on
$\Omega$ is denoted by ${\mcal M}(\Omega)$ while
 ${\mcal M}_0(\Omega)$ stands for the class of functions in ${\mcal M}(\Omega)$
that are finite $\mu_n$-a.e. on $\Omega$ and
${\mcal M}^+(\Omega)$ denotes the subset of ${\mcal M}(\Omega)$
consisting of all functions which are non-negative $\mu_n$-a.e. on $\Omega$.
When $\Omega$ is an interval $(a,b)\subseteq\mathbb{R}$, we denote
these sets by $\mathcal{M}(a,b)$,  $\mathcal{M}_0(a,b)$ and $\mathcal{M}^{+}(a,b),$
respectively. By $\mathcal{M}_0^{+}(a,b;\downarrow)$ we mean the
subset of $\mathcal{M}_0^{+}(a,b)$ containing all non-increasing
functions on the interval $(a,b)$ and  by $\mathcal{M}_0^{+}(a,b;\uparrow)$ we mean the
subset of $\mathcal{M}_0^{+}(a,b)$ containing all non-decreasing
functions on the interval $(a,b)$. The symbol ${\mcal W}(a,b)$
stands for the class of weight functions on $(a,b)\subseteq\bkR$
consisting of all $\mu_1$-measurable functions which are positive and finite
$\mu_1$-a.e. on $(a,b)$. The {\it non-increasing rearrangement} of  $f\in
{\mcal M}(\bkR^n)$ is the function $f^*$ defined by $
f^*(t):=\inf\left\{\lambda\geq 0:|\{x\in
\bkR^n:|f(x)|\!>\!\lambda\}|_n\leq t \right\}$ for all $t\geq 0$.
By
$f^{**}$ we denote the maximal function of $f^*$ given by
$f^{**}(t):=t^{-1}\int_0^tf^*(\tau)\,d\tau$, $t\in(0,+\infty)$. The maximal operator $f\mapsto f^{**}$ is subadditive (cf. \cite[p. 54]{BS:IO}). 


By $c$, $C$, $c_1$, $C_1$, $c_2$, $C_2$, etc.  we denote positive
constants independent of appropriate quantities. For two
non-negative expressions ({i.e.} functions or functionals) ${\mcal
A}$, ${\mcal B}$, the symbol ${\mcal A}\precsim {\mcal B}$ (or
${\mcal A}\succsim {\mcal B}$) means that $ {\mcal A}\leq c\, {\mcal
B}$ (or $c\,{\mcal A}\geq {\mcal B}$). If ${\mcal A}\precsim {\mcal
B}$ and ${\mcal A}\succsim{\mcal B}$, we write ${\mcal A}\approx
{\mcal B}$ and say that ${\mcal A}$ and ${\mcal B}$ are equivalent.
Throughout the paper we use the abbreviation $\text{LHS}(*)$
($\text{RHS}(*)$) for the left- (right-) hand side of relation
$(*)$. We adopt the convention that $a/(+\infty)=0$, 
$a/0=+\infty$ and $(+\infty)^a=+\infty$ for all $a\in(0,+\infty)$. If $p\in[1,+\infty]$, the conjugate
number $p'$ is given by $1/p+1/p'=1$. 
In the whole paper $\|.\|_{p;(c,d)},\, p\in
(0,+\infty]$, denotes the usual $L_p$-(quasi-)norm on the
interval $(c,d)\subseteq \bkR$.

We say that a positive, finite and
Lebesgue-measurable function $b$ is {\it slowly varying} on
$(0,+\infty)$, and write $b\in SV(0,+\infty)$, if, for each
$\varepsilon>0$, $t^{\varepsilon}b(t)$ is equivalent to a non-decreasing
function on $(0,+\infty)$ and  $t^{-\varepsilon}b(t)$ is equivalent to
a~non-increasing function on $(0,+\infty)$. Here we follow the definition of $SV(0,+\infty)$
given in \cite{GOT2002Lrrinwsvf}; for other definitions see, for example,
\cite{BGT87:RV}, \cite{EEv04:HOFSE}, \cite{EKP:OISRIQ}, and \cite{Nev02:LKSBRPE}. The
family $SV(0,+\infty)$ includes not only powers of iterated logarithms
and the broken logarithmic functions of \cite{EO00:RILFR} but also
such functions as $t\mapsto\exp\left(  \left\vert \log
t\right\vert ^{a}\right)  ,$ $a\in(0,1).$ (The last mentioned
function has the interesting property that it tends to infinity more
quickly than any positive power of the logarithmic function).

The replacement of the interval $(0,+\infty)$ in the definition of the class $SV(0,+\infty)$ by the interval $(0,1)$ yields the definition of the class $SV(0,1)$.

Let  $q\in(0,+\infty]$, $b\in
SV(0,+\infty)$,  $B_q(t):=\|\tau^{-1/q}b(\tau)\|_{q;(0,t)}$, $t\in(0,+\infty)$, and let $B_q(1)<+\infty$. Then, by
\cite[Lemma 2.1 (v)]{GO07:SEBtS},
\begin{equation*}\label{svcontinuousfunct0}
B_q\in SV(0,+\infty).
\end{equation*}

If $n\in\bkN$, $r\in(0,+\infty]$ and $b\in
SV(0,+\infty)$, then we put
\begin{equation*}\label{svcontinuousfunctnearinfty}
b_{r,n}(t):=\|\tau^{-1/r}b(\tau^{1/n})\|_{r;(t,+\infty)}\quad\text{for all}\quad t\in (0,+\infty).
\end{equation*}
Suppose $b_{r,n}(1)<+\infty$. By
\cite[Lemma 2.1 (ii), (v)]{GO07:SEBtS},
\begin{equation*}\label{svbrn}
b_{r,n}\in SV(0,+\infty).
\end{equation*}

By \cite[Lemma 2.2 (8)]{CGO2},
\begin{equation}\label{prop:SS_Limite_Frac_0}
\limsup_{x\rightarrow 0_+}\dfrac{\|\tau^{-1/q} b(\tau)\|_{q;(x,+\infty)}}{b(x)}=+\infty\quad\text{for any}\quad b\in SV(0,+\infty).
\end{equation}
By \cite[Lemma 2.1 (i))]{GO07:SEBtS} and the previous result, it also follows that,
\begin{equation}\label{prop:SS_Limite_Frac_Infty}
\limsup_{x\rightarrow +\infty}\dfrac{\|\tau^{-1/q} b(\tau)\|_{q;(0,x)}}{b(x)}=+\infty\quad\text{for any}\quad b\in SV(0,+\infty).
\end{equation}

More properties and examples of slowly varying functions can be
found in \cite[Chapter~V, p. 186]{Zyg57:TS}, \cite{BGT87:RV},
\cite{EKP:OISRIQ}, \cite{VojislavMaric:RVDE00}, \cite{Nev02:LKSBRPE},
\cite{GOT2002Lrrinwsvf} and \cite{GNO10:PotAnal}.

Let $n\in\bkN$, $p,q\in(0,+\infty]$, $b\in SV(0,+\infty)$.
The {\it Lorentz-Karamata} (LK)
space
 $L_{p,q;b}(\bkR^n)$ is
defined to be the set of all functions $f\in {\mcal M}(\bkR^n)$ such
that
\begin{equation}\label{eq:defquasinormLKspa}
\|f\|_{p,q;b}:=
\|t^{1/p-1/q}\;b(t)\;f^*(t)\|_{q;(0,+\infty)}<+\infty.
\end{equation}
 The local {\it Lorentz-Karamata}
space
 $L_{p,q;b}^{loc}(\bkR^n)$ consists of all functions $f\in {\mcal M}(\bkR^n)$ such
that
\begin{equation}\label{eq:defquasinormLKspaLocal}
\|f\|_{p,q;b}^{loc}:=
\|t^{1/p-1/q}\;b(t)\;f^*(t)\|_{q;(0,1)}<+\infty.
\end{equation}

Particular choices of $b$ give well-known spaces. If $m\in\bkN$,
$\balpha=(\alpha_1,\dots,\alpha_m)\in\bkR^m$ and
\[
b(t)=\bell^{\balpha}(t):=\prod_{i=1}^{m}\ell_{i}^{\alpha_{i}}(t)\text{
\ for all \ }t\in(0,+\infty)
\]
(where $ \ell(t)=\ell_{1}(t):=1+\left\vert \log t\right\vert ,\text{
}\ell_{i}(t):=\ell_{1}(\ell_{i-1} (t))\text{ if }i>1$), then the LK-space
$L_{p,q;b}(\bkR^n)$ is the generalized Lorentz-Zygmund space
$L_{p,q,\balpha}$ introduced in \cite{EGO4} and endowed with the
\mbox{(quasi-)}norm $\|f\|_{p,q;\balpha;\bkR^n}$,  which in turn becomes
the Lorentz-Zygmund space $L_{p,q}(\log L)^{\alpha_{1}}(\bkR^n)$ of Bennett
and Rudnick \cite{BR80:OLZS} when $m=1$. If $\balpha=(0,\dots, 0)$,
we obtain the Lorentz space $L_{p,q}(\bkR^n)$ endowed with the
(quasi-)norm $\|.\|_{p,q;\bkR^n}$, which is just the Lebesgue space
$L_p(\bkR^n)$ equipped with the {(quasi-)}norm
$\|.\|_{p}$ when $p=q$; if $p=q$ and $m=1$, we obtain the
Zygmund space $L_{p}(\log L)^{\alpha_1}(\bkR^n)$ endowed with the
(quasi-)norm $\|.\|_{p;\alpha_1}$.


Let $n\in\bkN$, $1< p< +\infty$, $1\leq q\leq +\infty$, $b\in SV(0,+\infty)$. 
The space $
Z_{p,q,n,b}(\bkR^n)$ is defined to be the set of all functions $f\in{\mcal M}(\bkR^n)$ such that
\begin{equation}\label{defquasinormZ:maintheo:1}
\|f\|_{Z_{p,q,n,b}}:=\big\|t^{-1/q}b(t^{1/n})\,\|f^*\|_{p;(0,t)}\big\|_{q;(0,+\infty)}<+\infty.
\end{equation}
We can define as well the local space $Z_{p,q,n,b}^{loc}(\bkR^n)$ as the set of all functions $f\in{\mcal M}(\bkR^n)$ for which 
\begin{equation}\label{defquasinormZ:maintheo:1Local}
\|f\|_{Z_{p,q,n,b}^{loc}}:=\big\|t^{-1/q}b(t^{1/n})\,\|f^*\|_{p;(0,t)}\big\|_{q;(0,1)}<+\infty.
\end{equation}

From \eqref{defquasinormZ:maintheo:1} and \eqref{defquasinormZ:maintheo:1Local} it is clear that
$$
Z_{p,q,n,b}(\bkR^n)\hookrightarrow Z_{p,q,n,b}^{loc}(\bkR^n).
$$

We refer to \cite{FK2004}, \cite{CF2005} and \cite{FFG2018} for the connection of these spaces with small Lebesgue spaces.

Let $h\in\bkR^n$. The first difference operator $\Delta_h$ is
defined on scalar functions $f$ on $\bkR^n$ by
$\Delta_hf(x)=f(x+h)-f(x)$ for all $x\in\bkR^n$.

Let  $p\in[1,+\infty)$. The {\it modulus of smoothness} of a function
$f$ in $L_p(\bkR^n)$ is given by
$$
\ds{\omega_1(f,t)_{p}:=\sup_{|h|\leq t}\|\Delta_hf\|_{p}\quad
\mbox{for all \ }t\geq 0}.
$$

\begin{definition}
Let $n\in\bkN$, $1\leq p<+\infty$, $1\leq r\leq +\infty$ and let $b\in SV(0,+\infty)$ be such that
\begin{equation}\label{cond:maintheo:1}
\|t^{-1/r}\;b(t)\|_{r;(0,1)}=+\infty\quad\text{and}\quad
\|t^{-1/r}\;b(t)\|_{r;(1,+\infty)}<+\infty.
\end{equation}
The Besov space $B_{p,r}^{0,b}(\bkR^n)$ consists of those $f\in L_p(\bkR^n)$ for which the norm
\begin{equation}\label{defi:normBesov}
\|f\|_{B_{p,r}^{0,b}}:=\|f\|_p+\|t^{-1/r}b(t)\omega_1(f,t)_p\|_{r;(0,+\infty)}
\end{equation}
is finite.

\end{definition}

\begin{remark}\label{remark:considerationsDefinBesov}
{\rm
\begin{trivlist}
\item[\hspace*{0.5cm}{\rm (i)}] An equivalent norm on $B_{p,r}^{0,b}(\bkR^n)$ is given by the functional
$$
\|f\|_{B_{p,r}^{0,b}}:=\|f\|_p+\|t^{-1/r}b(t)\omega_1(f,t)_p\|_{r;(0,1)}.
$$
We refer to \cite[{Remark 2.5 (iii)}]{CGO2} for more details.
\item[\hspace*{0.5cm}{\rm (ii)}] The assumption $ \|t^{-1/r}b(t)\|_{r;(1,+\infty)}<+\infty$ is natural. Otherwise the space $B_{p,r}^{0,b}(\bkR^n)$ is trivial (that is, it consists only of the zero element). We again refer to \cite[{Remark~2.5 (iii)}]{CGO2} for more details.

\item[\hspace*{0.5cm}{\rm (iii)}]  Note also that only the case when $\|t^{-1/r}b(t)\|_{r;(0,1)}=+\infty$ is of interest. Otherwise $B_{p,r}^{0,b}(\bkR^n)\equiv L_p(\bkR^n)$, cf. \cite[{Remark 2.5 (i)}]{CGO2}.

\item[\hspace*{0.5cm}{\rm (iv)}] An equivalent norm results on $B^{0,b}_{p,r}(\bkR^n)$ if the modulus of smoothness $\omega_1(f,\cdot)_p\,$ in \eqref{defi:normBesov} is replaced by the $k$-th order modulus of smoothness $\, \omega_k(f,\cdot)_p\,$, where $\, k \in \{ 2,3,4, \ldots \}$,  cf. \cite[{Remark 2.5 (ii)}]{CGO2}.
\end{trivlist}
}
\end{remark}

\begin{remark}\label{remark:13}
{\rm Assumption \eqref{cond:maintheo:1} implies that
\begin{equation}\label{eq:34}
\int_t^{+\infty}{\tau^{-1}b^r(\tau^{1/n})}\,d\tau\approx \int_t^{2}{\tau^{-1}b^r(\tau^{1/n})}\,d\tau\quad\text{for all}\quad  {t\in(0,1)}.
\end{equation}
Indeed, since
\begin{align*}
\int_t^{+\infty}{\tau^{-1}b^r(\tau^{1/n})}\,d\tau&<\int_t^{2}\tau^{-1}b^r(\tau^{1/n})\,d\tau+\int_1^{+\infty}\tau^{-1}b^r(\tau^{1/n})\,d\tau\\
&\precsim \int_t^{2}\tau^{-1}b^r(\tau^{1/n})\,d\tau+ 1\\
&\precsim \int_t^{2}\tau^{-1}b^r(\tau^{1/n})\,d\tau + \int_1^{2}\tau^{-1}b^r(\tau^{1/n})\,d\tau\\
&\precsim  \int_t^{2}\tau^{-1}b^r(\tau^{1/n})\,d\tau,\quad\text{for all}\quad  t\in(0,1),
\end{align*}
and since the reverse estimate is trivial, we see that \eqref{eq:34} holds.
}
\end{remark}

\section{Main Results}\label{sectionmain}

\begin{theorem}\label{theo:mainresult}
Let $n\in\bkN$, $1<p<+\infty$, $1\leq r, q<+\infty$, and let $b\in
SV(0,+\infty)$ be such that \eqref{cond:maintheo:1} holds.
If  \ $\bar{b}=\bar{b}_{(r,q)}\in SV(0,+\infty)$ is given by
\begin{equation}\label{deffunct:maintheo:1}
\bar{b}(t)=\bar{b}_{(r,q)}(t):=(b_r(t))^{1-r/q}(b(t))^{r/q}\quad\text{for all \ }t\in(0,+\infty),\mbox{\footnotemark$^)$}
\end{equation}
\footnotetext{$^{)}$ When $q=r$, $\bar{b}=\bar{b}_{(r,r)}=b$.}
where
\begin{equation}\label{deffunctbr:maintheo:1}
b_{r}(t):=\|\tau^{-1/r}b(\tau)\|_{r;(t,+\infty)}\quad\text{for all}\quad t\in (0,+\infty),
\end{equation}
then
\begin{equation}\label{embedding:maintheo:1}
B_{p,r}^{0,b}(\bkR^n)\hookrightarrow Z_{p,q,n,\bar{b}}(\bkR^n)=:Z(\bkR^n),
\end{equation}
if and only if $q\geq r$.

\end{theorem}

\begin{corollary}\label{cor:mainresult}
Let all the assumptions of \ {\rm Theorem \ref{theo:mainresult}} be satisfied and $r\leq q$. Let $\widetilde{b}=\widetilde{b}_{(r,q,n,p)}\in SV(0,+\infty)$ be given by
\begin{equation}\label{deffunct:maintcor:1}
\widetilde{b}(t)=\widetilde{b}_{(r,q,n,p)}(t):=\big[b_{r}(t^{1/n})\big]^{1-\frac{r}{q}+\frac{r}{\max\{p,q\}}} \big[b(t^{1/n})\big]^{\frac{r}{q}-\frac{r}{\max\{p,q\}}}\quad\text{for all \ }t\in(0,+\infty),
\end{equation}
with $b_r$ from \eqref{deffunctbr:maintheo:1}.
Then
\begin{equation}\label{embedding:maincor:1}
B_{p,r}^{0,b}(\bkR^n)\hookrightarrow L_{p,q;\widetilde{b}}(\bkR^n).
\end{equation}
\end{corollary}


\begin{remark}\label{remark:mainresult}\rm
\begin{trivlist}
\item[\hspace*{0.5cm}{\rm (i)}] Theorem \ref{theo:mainresult} gives in general a better result than Corollary \ref{cor:mainresult} since, if $1<p<+\infty$, $r\leq q$ and $q\neq p$, then the space $Z_{p,q,n,\bar{b}}(\bkR^n)$ is strictly smaller than the Lorentz-Karamata space $L_{p,q;\widetilde{b}}(\bkR^n)$ (cf. Remark \ref{remark:relationsZandL} below). By Theorem \ref{theorem24} below, the target spaces in \eqref{embedding:maintheo:1} and in \eqref{embedding:maincor:1} coincide if $1<p<+\infty$, $r\leq q$ and $q=p$.

\item[\hspace*{0.5cm}{\rm (ii)}] Embedding \eqref{embedding:maincor:1} implies that
\begin{equation}\label{embedding:mainremark:1}
B_{p,r}^{0,b}(\bkR^n)\hookrightarrow L_{p,q;\widetilde{b}}^{loc}(\bkR^n).
\end{equation}
The embedding \eqref{embedding:mainremark:1} was proved in \cite[{Theorem 3.3}]{CGO2} by a completely different method.  In \cite[Theorem 3.4]{CGO2} it was also shown that the target space $L_{p,r;\widetilde{b}}^{loc}(\bkR^n)$, {with $\widetilde{b}=\widetilde{b}_{(r,r,n,p)}$, given by \eqref{deffunct:maintcor:1} with $q$ replaced by $r$}, is optimal among the Lorentz-Karamata spaces $L_{P,Q;B}^{loc}(\bkR^n)$, with $1<P<+\infty$, $1\leq r\leq Q<+\infty$, $B\in SV(0,+\infty)$, for which the embedding \eqref{embedding:mainremark:1} holds.

%

\end{trivlist}
\end{remark}

\begin{theorem}\label{theo:relatedtolemma14}
Let all the assumptions of \ {\rm Theorem \ref{theo:mainresult}} be satisfied and $r\leq q$. Then the embedding
\begin{equation}\label{optimal:eq38:optimal}
B_{p,r}^{0,b}(\bkR^n)\hookrightarrow Z_{p,r,n,\bar{b}_{(r,r)}}(\bkR^n)=Z_{p,r,n,b}(\bkR^n)
\end{equation}
is optimal among all the embeddings
$$
B_{p,r}^{0,b}(\bkR^n)\hookrightarrow Z_{p,q,n,\bar{b}_{(r,q)}}(\bkR^n)\quad\text{with}\quad r\leq q<+\infty.
$$
\end{theorem}

Concerning the sharpness of embedding \eqref{embedding:maintheo:1} with respect to the function $\bar{b}$, we have the following result.

\begin{theorem}\label{sharpnesswithrespecttofunction}
Let all the assumptions of \ {\rm Theorem \ref{theo:mainresult}} be satisfied and $r\leq q$. Let $\kappa\in{{\mcal M}_0^+(0,+\infty;\downarrow)} \cap SV(0,+\infty)$ be such that
\begin{equation}\label{eq:41}
B_{p,r}^{0,b}(\bkR^n)\hookrightarrow Z_{p,q,n,\kappa\bar{b}}(\bkR^n),
\end{equation}
with $\bar{b}$ given by \eqref{deffunct:maintheo:1}.
Then the function $\kappa$ is bounded on the interval $(0,+\infty)$.
\end{theorem}

\begin{remark}
{\rm If all the assumptions of \ {\rm Theorem \ref{theo:mainresult}} are satisfied and $r\leq q$, then embedding \eqref{eq:41} can be replaced by
\begin{equation}\label{eq:41Local}
B_{p,r}^{0,b}(\bkR^n)\hookrightarrow Z_{p,q,n,\kappa\bar{b}}^{loc}(\bkR^n).
\end{equation}
Indeed, \eqref{eq:41} implies \eqref{eq:41Local}. On the other hand, \eqref{eq:41Local} means that
\begin{equation}\label{estimate:localversusGlobal0}
\big\|t^{-1/q}\kappa(t)\,\bar{b}(t^{1/n})\,\|f^*\|_{p;(0,t)}\big\|_{q;(0,1)}\precsim \|f\|_{B_{p,r}^{0,b}}\quad\text{for all}\quad f\in B_{p,r}^{0,b}(\bkR^n).
\end{equation}
Moreover, since $\kappa$ is non-increasing, \eqref{embedding:maintheo:1} also implies that
\begin{align}
\big\|t^{-1/q}\kappa(t)\,\bar{b}(t^{1/n})\,\|f^*\|_{p;(0,t)}\big\|_{q;(1,+\infty)}&\leq \kappa(1)\, \big\|t^{-1/q}\bar{b}(t^{1/n})\,\|f^*\|_{p;(0,t)}\big\|_{q;(1,+\infty)}\label{estimate:localversusGlobal}\\
&\leq \kappa(1)\, \big\|t^{-1/q}\bar{b}(t^{1/n})\,\|f^*\|_{p;(0,t)}\big\|_{q;(0,+\infty)}\nonumber\\
&\precsim \|f\|_{B_{p,r}^{0,b}}\quad\text{for all}\quad f\in B_{p,r}^{0,b}(\bkR^n).\nonumber
\end{align}
Consequently, \eqref{estimate:localversusGlobal0} and \eqref{estimate:localversusGlobal} imply \eqref{eq:41}.
}
\end{remark}

The next result show that if $q\neq p$, then the target space $Z_{p,q,n,\bar{b}}(\bkR^n)$ in \eqref{embedding:maintheo:1} is strictly smaller than the target space $L_{p,q;\widetilde{b}}(\bkR^n)$ in \eqref{embedding:maincor:1}.


\begin{theorem}\label{theorem24}
Let all the assumptions of \ {\rm of Corollary \ref{cor:mainresult}} be satisfied and $r\leq q$. Assume that $\bar{b},\widetilde{b}\in SV(0,+\infty)$ are given by \eqref{deffunct:maintheo:1} and \eqref{deffunct:maintcor:1}, respectively. If the space $Z_{p,q,n,\bar{b}}(\bkR^n)$ is defined as in {\rm Theorem \ref{theo:mainresult}}, then
\begin{equation}\label{eq:61aux1}
Z_{p,q,n,\bar{b}}(\bkR^n)\hookrightarrow L_{p,q;\widetilde{b}}(\bkR^n).
\end{equation}
Moreover,
\begin{equation}\label{eq:61aux3}
Z_{p,q,n,\bar{b}}(\bkR^n)=L_{p,q;\widetilde{b}}(\bkR^n)\quad\text{if}\quad q=p,
\end{equation}
and
\begin{equation}\label{eq:62}
 L_{p,q;\widetilde{b}}(\bkR^n)\not\hookrightarrow Z_{p,q,n,\bar{b}}(\bkR^n)\quad\text{if}\quad q\neq p.
\end{equation}
\end{theorem}

\begin{remark}\label{remark:relationsZandL}
{\rm Let all the assumptions of Corollary \ref{cor:mainresult}  be satisfied. Since then $Z_{p,q,n,\bar{b}}(\bkR^n)$ and $L_{p,q;\widetilde{b}}(\bkR^n)$ are r.i. Banach function spaces,
\cite[Chap. I, p.7, Theorem 1.8]{BS:IO} and \eqref{eq:62} imply that
 $$
 L_{p,q;\widetilde{b}}(\bkR^n)\not\subset Z_{p,q,n,\bar{b}}(\bkR^n)\quad\text{if}\quad q\neq p.
 $$
This, together with \eqref{eq:61aux1}, gives
\begin{equation}\label{eq:61}
Z_{p,q,n,\bar{b}}(\bkR^n)\stackrel[\neq]{}{\hookrightarrow}  L_{p,q;\widetilde{b}}(\bkR^n)\quad\text{if}\quad q\neq p.
\end{equation}
 }
\end{remark}

\section{Auxiliary assertions}

\begin{lemma}\label{Lemma10:aux}
If \ $1\leq p<+\infty$, $1\leq r\leq q <+\infty$, $b\in SV(0,+\infty)$ and  $\bar{b}\in SV(0,+\infty)$ is given by
\begin{equation*}\label{deffunct:lemmaaux:1}
\bar{b}(t):=(b_r(t))^{1-r/q}(b(t))^{r/q}\quad\text{for all}\quad t\in(0,+\infty),
\end{equation*}
then
\begin{equation}\label{embedding:lemmaaux:1}
B_{p,r}^{0,b}(\bkR^n)\hookrightarrow B_{p,q}^{0,\bar{b}}(\bkR^n).
\end{equation}
\end{lemma}
\bpr Emebdding \eqref{embedding:lemmaaux:1} means that, for all $f\in B_{p,r}^{0,b}(\bkR^n)$,
$$
\|f\|_p+\|t^{-1/q}\bar{b}(t)\omega_1(f,t)_p\|_{q;(0,+\infty)}\precsim \|f\|_p+\|t^{-1/r}b(t)\omega_1(f,t)_p\|_{r;(0,+\infty)}.
$$
Thus, \eqref{embedding:lemmaaux:1} holds if,  for all $f\in B_{p,r}^{0,b}(\bkR^n)$, 
\begin{equation}\label{proofineq:lemmaaux:1}
\|t^{-1/q}\bar{b}(t)\omega_1(f,t)_p\|_{q;(0,+\infty)}\precsim \|t^{-1/r}b(t)\omega_1(f,t)_p\|_{r;(0,+\infty)}.
\end{equation}
Since $\omega_1(f,\cdot)_p\in\mathcal{M}_0^+(0,+\infty;\uparrow)$ for any $f\in B_{p,r}^{0,b}(\bkR^n)$, estimate \eqref{proofineq:lemmaaux:1} will be satisfied provided that
\begin{equation}\label{proofineq13:lemmaaux:1}
\left(\int_0^{+\infty}w(t)g^q(t)\,dt\right)^{1/q}\precsim \left(\int_0^{+\infty}v(t)g^r(t)\,dt\right)^{1/r}\quad\text{for all}\quad g\in \mathcal{M}_0^+(0,+\infty;\uparrow),
\end{equation}
where
\begin{equation}\label{proofineq14:lemmaaux:1}
w(t):=t^{-1}{\bar{b}}^q(t)\quad\text{and}\quad v(t):=t^{-1}b^r(t)\quad\text{for all}\quad t\in(0,+\infty).
\end{equation}
Putting
$$
W(t):=\int_t^{+\infty}w(\tau)\,d\tau\quad\text{and}\quad V(t):=\int_t^{+\infty}v(\tau)\,d\tau\quad\text{for all}\quad t\in(0,+\infty),
$$
we obtain by \cite[Proposition 2.1]{H-Step} that inequality \eqref{proofineq13:lemmaaux:1} holds if
\begin{equation}\label{proofineq15:lemmaaux:1}
\sup_{t\in(0,+\infty)}W^{1/q}(t)V^{-1/r}(t)<+\infty.
\end{equation}
Since condition \eqref{proofineq15:lemmaaux:1} is satisfied when $w$ and $v$ are given by  \eqref{proofineq14:lemmaaux:1}, the result follows.
\epr

\bigskip

We shall also need the result mentioned in \cite[Chap. V, Corollary 4.20, p. 346]{BS:IO}, which states that
\begin{equation}\label{eq:17:embBintoLinfty}
B_{p,1}^{n/p}(\bkR^n)\hookrightarrow L_{\infty}(\bkR^n)\quad\text{if $1\leq p<+\infty$ and $n\in\bkN$}.
\end{equation}

\medskip

\begin{definition}\label{defin:1:interpolatonspaces}
Let $(X_0,X_1)$ be a compatible couple of quasi-Banach spaces. 
\begin{trivlist}
\item[\hspace*{0.5cm}{\rm (i)}] For each $t\in(0,+\infty)$, $K(\cdot,t;X_0,X_1)$ is the Peetre's $K$-functional defined by
$$
K(f,t;X_0,X_1):=\inf_{f=f_0+f_1}(\|f_0\|_{X_0}+t\|f_1\|_{X_1})\quad\text{for any $f\in X_0+X_1$}.
$$
Sometimes, we denote $K(f,t;X_0,X_1)$ simply by $K(f,t)$.

\item[\hspace*{0.5cm}{\rm (ii)}] For $0\leq \theta\leq 1$, $0<q\leq +\infty$, and $b\in SV(0,+\infty)$, we put
\begin{equation}\label{eq:circ_1}
\overline{X}_{\theta,q,b}\equiv(X_0,X_1)_{\theta,q,b}:=\{f\in X_0+X_1: \|f\|_{\theta,q,b}<+\infty\},
\end{equation}
where
\begin{equation}\label{eq:circ_2}
 \|f\|_{\theta,q,b}:=\|t^{-\theta-1/q}\,b(t)\,K(f,t;X_0,X_1)\|_{q;(0,+\infty)}.
\end{equation}
If $b\equiv 1$ on $(0,+\infty)$, we write $(X_0,X_1)_{\theta,q}$ instead of $(X_0,X_1)_{\theta,q,b}$.

\item[\hspace*{0.5cm}{\rm (iii)}]  For $0<\sigma< 1$, $0<r,q\leq +\infty$, and $a, b\in SV(0,+\infty)$, we define the space
\begin{equation}\label{eq:circ_3}
\overline{X}_{\sigma,r,b,q,a}^{\mcal L}\equiv(X_0,X_1)_{\sigma,r,b,q,a}^{\mcal L}:=\{f\in X_0+X_1: \|f\|_{{\mcal L};\sigma,r,b,q,a}<+\infty\},
\end{equation}
where
\begin{equation}\label{eq:circ_4}
 \|f\|_{{\mcal L};\sigma,r,b,q,a}:=\left\|t^{-1/r}\,\frac{b(t)}{a(t)}\|\tau^{-\sigma-1/q}\,a(\tau)\,K(f,\tau;X_0,X_1)\|_{q;(0,t)}\right\|_{r;(0,+\infty)}.
\end{equation}
{\rm (}The superscript and the subscript ${\mcal L}$ in \eqref{eq:circ_3} and \eqref{eq:circ_4}, respectively, is an indication of the fact that the local {\rm (}quasi-{\rm )}norm in  \eqref{eq:circ_4} is taken from the left end of the interval $(0,+\infty)$.{\rm )}
\end{trivlist}
\end{definition}

Next we collect some weighted inequalities, which are needed in the rest of the paper. 
To prove Theorem \ref{Theorem5:circ} below (involving the embedding $\overline{X}_{\theta,r,b,q,a}^{\mcal L}\hookrightarrow \overline{X}_{\theta,r,d}$), we shall make use of the following assertions.

\begin{theorem}[{\cite[Theorem 2.1]{Lai:1993wt}}]\label{theorem3=Theorem2.1:Lai:1993wt}
Let $1\leq P\leq Q<+\infty$, $\varphi:\bkR_+\times\bkR_+\longrightarrow \bkR_+$, and $v, w\in {\mcal W}(0,+\infty)$. Then the inequality
\begin{equation}\label{eq7:circ}
\left(\int_0^{+\infty}h^Q(x)\,v(x)\,dx\right)^{1/Q}\leq C \left(\int_0^{+\infty}\left(\int_0^{+\infty} \varphi(x,y)\,h(y)\,dy\right)^P\,w(x)\,dx\right)^{1/P}
\end{equation}
holds for all $h\in{\mcal M}_0^+(0,+\infty;\downarrow)$ if and only if, for all $z\in(0,+\infty)$,
\begin{equation}\label{eq8:circ}
\left(\int_0^{z}v(x)\,dx\right)^{1/Q}\leq C \left(\int_0^{+\infty}\left(\int_0^{z} \varphi(x,y) \,dy\right)^P\,w(x)\,dx\right)^{1/P}.
\end{equation}
{\rm (}The constant C in both inequalities is the same.{\rm )}
\end{theorem}

\begin{lemma}[{\cite[Lemma 3.4]{EO00:RILFR}}]\label{lemma4=Lemma 3.4:EO00:RILFR}
Let $s\in(0,1)$, $\varphi, w\in{\mcal W}(0,+\infty)$ and define $v\in{\mcal W}(0,+\infty)$ by
\begin{equation}\label{eq9:circ}
v(y)=w^{1-s}(y)\left(\varphi(y)\int_y^{+\infty}w(x)\,dx\right)^s,\quad y\in(0,+\infty).
\end{equation}
Then, for all $h\in{\mcal M}_0^+(0,+\infty)$,
\begin{equation}\label{eq10:circ}
\int_0^{+\infty}h^s(x)\,v(x)\,dx\leq s^s \int_0^{+\infty}\left(\int_0^{x} \varphi(y) h(y)\,dy\right)^s\,w(x)\,dx.
\end{equation}
\end{lemma}

\

To prove Theorem \ref{Theorem5:circ:STAR} mentioned below (which is a counterpart of Theorem \ref{Theorem5:circ}), we shall make use of the following assertions.

\begin{theorem}[{\cite[Theorem 2.2]{Lai:1993wt}}]\label{theorem3=Theorem2.1:Lai:1993wt:STAR}
Let $0<Q\leq P\leq 1$, $\varphi:\bkR_+\times\bkR_+\longrightarrow \bkR_+$, and $v, w\in {\mcal W}(0,+\infty)$. Then the inequality
\begin{equation}\label{eq7:circ:STAR}
\left(\int_0^{+\infty}\left(\int_0^{+\infty} \varphi(x,y) \,h(y)\,dy\right)^P\,w(x)\,dx\right)^{1/P}\leq C \left(\int_0^{+\infty}h^Q(x)\,v(x)\,dx\right)^{1/Q}
\end{equation}
holds for all $h\in{\mcal M}_0^+(0,+\infty;\downarrow)$ if and only if, for all $z\in(0,+\infty)$,
\begin{equation}\label{eq8:circ:STAR}
 \left(\int_0^{+\infty}\left(\int_0^{z} \varphi(x,y) \,dy\right)^P\,w(x)\,dx\right)^{1/P}\leq C \left(\int_0^{z}v(x)\,dx\right)^{1/Q}.
\end{equation}
{\rm (}The constant C in both inequalities is the same.{\rm )}
\end{theorem}

\begin{theorem}[{\cite[Theorem 5.9, p. 63]{OK90:HTI}}]\label{Theorem20=Theorem5.9:OK90:HTI}
Let $1\leq P\leq Q\leq +\infty$ and $v, w\in {\mcal W}(0,+\infty)$.
Then, the inequality
\begin{equation}\label{eq:50}
\left\|\left(\int_0^xg(y)\,dy\right)\,w^{1/Q}(x)\right\|_{Q;(0,+\infty)}\precsim \|g(x)\,v^{1/P}(x)\|_{P;(0,+\infty)}
\end{equation}
holds for all $g\in{\mcal M}_0^+(0,+\infty)$ if and only if
\begin{equation}\label{eq:51}
\sup_{x\in(0,+\infty)}\|w^{1/Q}\|_{Q;(x,+\infty)}\,{\|v^{-1/P}\|_{P';(0,x)}}<+\infty.
 \end{equation}
\end{theorem}

To study the relations among the target spaces in 
Theorem~\ref{theo:mainresult} and in Corollary~\ref{cor:mainresult}, we shall need the following result.

\begin{theorem}[{\cite[Theorem 2]{Swa90:BCOCLS}}]\label{theorem23=theorem2Sawyer}
Suppose $1<P\leq Q<+\infty$ and $v,w\in{\mcal W}(0,+\infty)$. Then the inequality
\begin{equation}\label{eq:58}
\left(\int_0^{+\infty}\left(x^{-1}\int_0^xg(t)\,dt\right)^Qw(x)\,dx\right)^{1/Q}\leq C 
\left(\int_0^{+\infty} g^P(x)\,v(x)\,dx\right)^{1/P}
\end{equation}
holds for all $g\in{\mcal M}_0^+(0,+\infty;\downarrow)$ if and only if there are constants $A$ and $B$ such that
\begin{equation}\label{eq:59}
\left(\int_0^{z}w(x)\,dx\right)^{1/Q}\leq A 
\left(\int_0^{z} v(x)\,dx\right)^{1/P}
\end{equation}
and
\begin{equation}\label{eq:60}
\left(\int_z^{+\infty}x^{-Q}\,w(x)\,dx\right)^{1/Q}
\left(\int_0^{z}\left(\frac{x}{\int_0^{x}v(t)\,dt}\right)^{P'}\,v(x)\,dx\right)^{1/P'}\leq B
\end{equation}
for all $z\in(0,+\infty)$.

Moreover, if $C$ is the best constant in \eqref{eq:58}, then $C\approx A+B$.
\end{theorem}

The next result extends \cite[Theorem 4.7]{EO00:RILFR}.

\begin{theorem}\label{Theorem5:circ}
Let $(X_0,X_1)$ be a compatible couple of quasi-Banach spaces. If \ $\theta\in [0,1)$,  $0<r,q<+\infty$ and $a,b\in SV (0,+\infty)$, then 
\begin{equation}\label{eq:11:circ}
\overline{X}_{\theta,r,b,q,a}^{\mcal L}\hookrightarrow \overline{X}_{\theta,r,d}:=(X_0,X_1)_{\theta,r,d},
\end{equation}
where
\begin{equation}\label{eq:12:circ}
d(x):= b(x)\left(\frac{\int_x^{+\infty}y^{-1}\left(\frac{b(y)}{a(y)}\right)^r\,dy}{\left(\frac{b(x)}{a(x)}\right)^r}\right)^{\frac{1}{\max\{q,r\}}},\quad x\in(0,+\infty).
\end{equation}
\end{theorem}
\bpr
Embedding \eqref{eq:11:circ} means that, for all $f\in \overline{X}_{\theta,r,b,q,a}^{\mcal L}$,
\begin{equation}\label{eq:13:circ}
\|t^{-\theta-1/r}\,d(t)\,K(f,t)\|_{r;(0,+\infty)}\precsim  
\left\|t^{-1/r}\,\frac{b(t)}{a(t)}\|\tau^{-\theta-1/q}\,a(\tau)\,K(f,\tau)\|_{q;(0,t)}\right\|_{r;(0,+\infty)}.
\end{equation}

Put $Q:=\frac{r}{q}$ and $h(t):=[K(f,t)/t]^q$, $t\in (0,+\infty)$. Then $h\in{\mcal M}_0^+(0,+\infty;\downarrow)$ and inequality
\eqref{eq:13:circ} can be rewritten  as
$$
\left(\int_0^{+\infty}\!\!t^{(1-\theta)r-1}\,d^r(t)\,h^Q(t)\,dt\right)^{1/Q}\!\!\!\precsim  \!
\left(\int_0^{+\infty}\!\!t^{-1}\!\left(\frac{b(t)}{a(t)}\right)^r\!\!\left(\int_0^{t}\!\!\tau^{(1-\theta)q-1}\,a^q(\tau)\,h(\tau)\,d\tau\!\right)^Q\!\!dt\right)^{1/Q}\!\!,
$$
that is
\begin{equation}\label{eq:15:circ}
\left(\int_0^{+\infty}v(x)\,h^Q(x)\,dx\right)^{1/Q}\precsim  
\left(\int_0^{+\infty}w(x)\left(\int_0^{x}y^{(1-\theta)q-1}\,a^q(y)\,h(y)\,dy\right)^Q\,dx\right)^{1/Q},
\end{equation}
where
$v(x):=x^{(1-\theta)r-1}\,d^r(x)$, $w(x):=x^{-1}\,\left(\frac{b(x)}{a(x)}\right)^r$, $x\in(0,+\infty)$, and $h\in{\mcal M}_0^+(0,+\infty;\downarrow)$.

\begin{trivlist}
\item[\hspace*{0.5cm}{\rm I.}] Assume first that $Q=\frac{r}{q}\in[1,+\infty)$. Put $P:=Q$ and
$$
\varphi(x,y):=\chi_{(0,x)}(y)y^{(1-\theta)q-1}\,a^q(y),\quad x,y\in(0,+\infty).
$$
Then inequality \eqref{eq:15:circ} coincides with inequality \eqref{eq7:circ}. Since one can verify that \eqref{eq8:circ} is satisfied in our case, Theorem \ref{theorem3=Theorem2.1:Lai:1993wt} implies that inequality  \eqref{eq:15:circ} holds for all $h\in{\mcal M}_0^+(0,+\infty;\downarrow)$. Thus, embedding \eqref{eq:11:circ} is satisfied if $0<q\leq r<+\infty$.

\item[\hspace*{0.5cm}{\rm II.}] Assume now that $Q=\frac{r}{q}\in(0,1)$. Put $s:=Q$ and $
\varphi(x):=x^{(1-\theta)q-1}\,a^q(x)$, $x\in(0,+\infty)$. Then inequality \eqref{eq:15:circ} is of the same form as inequality \eqref{eq10:circ}. Since in our case one can verify that
$$
v(y)\approx w^{1-s}(y)\left(\varphi(y)\int_y^{+\infty}w(x)\,dx\right)^s\quad\text{for all}\quad y\in(0,+\infty),
$$
Lemma \ref{lemma4=Lemma 3.4:EO00:RILFR} implies that inequality \eqref{eq:15:circ} is satisfied for all $h\in{\mcal M}_0^+(0,+\infty;\downarrow)$. Consequently, embedding \eqref{eq:11:circ} holds if $0<r< q<+\infty$.
\end{trivlist}
\epr

\begin{remark}{\rm
If we assume that the compatible couple $(X_0,X_1)$ of (quasi-) Banach spaces is ordered in the sense that $X_1\hookrightarrow X_0$, then one can define the spaces $(X_0,X_1)_{\theta,q,b}$ and $(X_0,X_1)_{\theta,r,b,q,a}^{\mcal L}$ just as in Definition \ref{defin:1:interpolatonspaces}, except that the role of the interval $(0,+\infty)$ is played by the interval $(0,1)$.

Moreover, one can show that when an ordered couple $(X_0,X_1)$ is considered, Theorem~\ref{Theorem5:circ} remain true provided that the functions $a,b\in SV(0,1)$ are extended to the interval $(0,2)$ by setting $a(x){\approx} b(x){\approx }1$, $x\in[1,2)$, and the function $d$ is given by
\begin{equation}\label{eq:12:circ:star}
d(x):= b(x)\left(\frac{\int_x^{2}y^{-1}\left(\frac{b(y)}{a(y)}\right)^r\,dy}{\left(\frac{b(x)}{a(x)}\right)^r}\right)^{\frac{1}{\max\{q,r\}}},\quad x\in(0,1).
\end{equation}
}
\end{remark}

Now we are going to prove a counterpart of Theorem \ref{Theorem5:circ}, which extends \cite[Lemma~4.3]{EO00:RILFR}.

\begin{theorem}\label{Theorem5:circ:STAR}
Let $(X_0,X_1)$ be a compatible couple of quasi-Banach spaces. If  \ $\theta\in [0,1)$, $0<r,q<+\infty$ and $a,b\in SV (0,+\infty)$, then 
\begin{equation}\label{eq:11:circ:STAR}
 \overline{X}_{\theta,r,d}:=(X_0,X_1)_{\theta,r,d}\hookrightarrow\overline{X}_{\theta,r,b,q,a}^{\mcal L},
\end{equation}
where
\begin{equation}\label{eq:12:circ:STAR}
d(x):= b(x)\left(\frac{\int_x^{+\infty}y^{-1}\left(\frac{b(y)}{a(y)}\right)^r\,dy}{\left(\frac{b(x)}{a(x)}\right)^r}\right)^{\frac{1}{\min\{q,r\}}},\quad x\in(0,+\infty).
\end{equation}
\end{theorem}
\bpr
Embedding \eqref{eq:11:circ:STAR} means that, for all $f\in \overline{X}_{\theta,r,d}$,
\begin{equation}\label{eq:53}
\left\|t^{-1/r}\,\frac{b(t)}{a(t)}\|\tau^{-\theta-1/q}\,a(\tau)\,K(f,\tau)\|_{q;(0,t)}\right\|_{r;(0,+\infty)}\precsim \|t^{-\theta-1/r}\,d(t)\,K(f,t)\|_{r;(0,+\infty)}.
\end{equation}

Put $Q:=\frac{r}{q}$ and $h(t):=[K(f,t)/t]^q$, $t\in (0,+\infty)$. Then $h\in{\mcal M}_0^+(0,+\infty;\downarrow)$ and inequality
\eqref{eq:53} can be rewritten  as
\begin{align}\label{eq:54}
\left(\int_0^{+\infty}t^{-1}\,\left(\frac{b(t)}{a(t)}\right)^r\right.&\left.\left(\int_0^{t}\tau^{(1-\theta)q-1}\,a^q(\tau)\,h(\tau)\,d\tau\right)^Q\,dt\right)^{1/Q}\\
&\precsim \left(\int_0^{+\infty}t^{(1-\theta)r-1}\,d^r(t)\,h^Q(t)\,dt\right)^{1/Q}.\nonumber
\end{align}

\begin{trivlist}
\item[\hspace*{0.5cm}{\rm I.}] Assume first that $Q:=\frac{r}{q}\in (0,1]$. Remark that  \eqref{eq:54} can be formulated as
\begin{equation}\label{eq:55}
\left(\int_0^{+\infty}w(x)\left(\int_0^{x}y^{(1-\theta)q-1}\,a^q(y)\,h(y)\,dy\right)^Q\,dx\right)^{1/Q}\precsim  \left(\int_0^{+\infty}v(x)\,h^Q(x)\,dx\right)^{1/Q},
\end{equation}
where
$v(x):=x^{(1-\theta)r-1}\,d^r(x)$, $w(x):=x^{-1}\,\left(\frac{b(x)}{a(x)}\right)^r$, $x\in(0,+\infty)$, and $h\in{\mcal M}_0^+(0,+\infty;\downarrow)$.
Put $P:=Q$ and
$$
\varphi(x,y):=\chi_{(0,x)}(y)y^{(1-\theta)q-1}\,a^q(y),\quad x,y\in(0,+\infty).
$$
Then inequality \eqref{eq:55} coincides with inequality \eqref{eq7:circ:STAR}. Using properties of slowly varying functions, one can verify that \eqref{eq8:circ:STAR} is satisfied in our case. Thus, Theorem \ref{theorem3=Theorem2.1:Lai:1993wt:STAR} implies that inequality  \eqref{eq:55} holds for all $h\in{\mcal M}_0^+(0,+\infty;\downarrow)$. Consequently, embedding \eqref{eq:11:circ:STAR} is satisfied if $0<r\leq q<+\infty$.

\item[\hspace*{0.5cm}{\rm II.}] Suppose now that $Q=\frac{r}{q}\in(1,+\infty)$. Putting $
g(\tau):=\tau^{(1-\theta)q-1}\,a^q(\tau)h(\tau)$, $\tau\in(0,+\infty)$, we have $g\in {\mcal M}_0^+(0,+\infty)$ and inequality \eqref{eq:54} can be rewritten as  
$$
\left(\int_0^{+\infty}t^{-1}\left(\frac{b(t)}{a(t)}\right)^r\left(\int_0^{t}g(\tau)\,d\tau\right)^Q\,dt\right)^{1/Q}\precsim  \left(\int_0^{+\infty}
t^{\frac{r}{q}-1}\left(\frac{d(t)}{a(t)}\right)^r\,g^Q(t)\,dt\right)^{1/Q},
$$
i.e., as \eqref{eq:50}, with $P:=Q$,
$v(t):=t^{\frac{r}{q}-1}\left(\frac{d(t)}{a(t)}\right)^r$,  $w(t):=t^{-1}\left(\frac{b(t)}{a(t)}\right)^r$, $t\in(0,+\infty)$.
One can verify that, for all $x\in(0,+\infty)$,
$$
\|w^{1/Q}\|_{Q;(x,+\infty)}=\left(\int_x^{+\infty}t^{-1}\left(\frac{b(t)}{a(t)}\right)^r\,dt\right)^{q/r}
$$
and
$$
 \|v^{-1/P}\|_{P';(0,x)}= \|v^{-1/Q}\|_{Q';(0,x)}\precsim \left(\int_x^{+\infty}t^{-1}\left(\frac{b(t)}{a(t)}\right)^r\,dt\right)^{-q/r}.
$$
Consequently, condition \eqref{eq:51} is satisfied, which means that inequality \eqref{eq:54} holds for all $g\in {\mcal M}_0^+(0,+\infty)$. Hence, embedding \eqref{eq:11:circ:STAR} holds as well if $0<q<r<+\infty$.
\end{trivlist}
\epr

\begin{remark}
If \ $r=q$, then {\rm Theorems  \ref{Theorem5:circ}} and {\rm \ref{Theorem5:circ:STAR}} imply that
\begin{equation}\label{eq:56}
 \overline{X}_{\theta,r,d}=\overline{X}_{\theta,r,b,r,a}^{\mcal L},
\end{equation}
with $d$ given by \eqref{eq:12:circ} {\rm (}or by  \eqref{eq:12:circ:STAR}{\rm )}, for all $\theta\in [0,1)$, $0<r<+\infty$ and $a,b\in SV(0,+\infty)$.
\end{remark}

\begin{theorem}[{\cite[Chap. V, Corollary 4.13, p. 341]{BS:IO}}]
Suppose $1\leq p\leq +\infty$ and $r\in\bkN$. If \ $0<\theta<1$ and $1\leq q\leq +\infty$, then
\begin{equation}\label{eq:16:interpLpW=B}
(L_p(\bkR^n),W_p^{r}(\bkR^n))_{\theta,q}=B_{p,q}^{\theta r}(\bkR^n),
\end{equation}
with equivalent norms.
\end{theorem}

We refer to \cite[p. 310] {BS:IO}, in the Banach setting, and to \cite[Definition 3.12]{EOP02:RILF}, in the general setting, for the definition of intermediate spaces of $(X_0, X_1)$ of class 0 or class~1.
 
\begin{theorem}[{\cite[Theorem 3.5]{GOT2002Lrrinwsvf}}]\label{Theo3.5:GOT2002Lrrinwsvf}
Let $0<\theta_0<\theta_1<1$, $0<q,q_0,q_1\leq +\infty$ and let $b,b_0,b_1\in SV(0,+\infty)$. Put $\overline{X}_{\theta_i}:=(X_0,X_1)_{\theta_i, q_i, b_i}$ for $i\in\{0,1\}$. Suppose that $\overline{X}_i$, $i\in\{0,1\}$, are intermediate spaces between $X_0$ and $X_1$ of class $i$. Then
\begin{align}
(\overline{X}_0,\overline{X}_{\theta_1})_{0,q,b}&=(X_0,X_1)_{0,q,b\circ \rho},\quad\text{where}\quad \rho(t)=t^{\theta_1}/ {b_1(t)},\quad t\in(0,+\infty),\label{eq:3.20:GOT2002}\\
(\overline{X}_{\theta_0},\overline{X}_1)_{0,q,b}&=\overline{X}_{\theta_0,q,b_0 b\circ\rho, q_0,b_0}^{\mcal L},\quad\text{where}\quad \rho(t)=t^{1-\theta_0}{b_0(t)},\quad t\in(0,+\infty). \label{eq:3.21:GOT2002}
\end{align}
\end{theorem}

\begin{lemma}[{\cite[Lemma 5.2]{GOT2002Lrrinwsvf}}]\label{lemma5.2:GOT2002Lrrinwsvf}
Let $0<\theta<1$, $0<q \leq +\infty$, $0<s<+\infty$, and let $b\in SV(0,+\infty)$. Then, for all $f\in L_s+L_{\infty}$ and all $t\in (0,+\infty]$,
\begin{equation}\label{eq:5.10:GOT2002}
\|\tau^{-\theta-1/q}\, b(\tau)\,K(f,\tau;L_s,L_{\infty})\|_{q;(0,t)}\approx\|y^{(1-\theta)/s-1/q}\,b(y^{1/s})\,f^*(y)\|_{q;(0,t^s)},
\end{equation}
and, for all $f\in L_{s,\infty}+L_{\infty}$ and all $t\in (0,+\infty]$,
\begin{equation}\label{eq:5.11:GOT2002}
\|\tau^{-\theta-1/q}\, b(\tau)\,K(f,\tau;L_{s,\infty},L_{\infty})\|_{q;(0,t)}\approx\|y^{(1-\theta)/s-1/q}\,b(y^{1/s})\,f^*(y)\|_{q;(0,t^s)}.
\end{equation}
\end{lemma}

\section{Proofs of Theorem \ref{theo:mainresult} and Corollary \ref{cor:mainresult}}
To prove the sufficient part of our first main result (i.e. Theorem \ref{theo:mainresult}), we shall use the following assertion.

\begin{theorem}\label{theorem11}
Let $n\in\bkN$, $1<p<+\infty$, $1\leq q<+\infty$, and let $\bar{b}\in SV(0,+\infty)$ be such that
\begin{equation}\label{eq:18}
\|t^{-1/q}\bar{b}(t)\|_{q;(0,1)}=+\infty\quad \text{and}\quad \|t^{-1/q}\bar{b}(t)\|_{q;(1,+\infty)}<+\infty.
\end{equation}
Then
\begin{equation}\label{emb:19}
B_{p.q}^{0,\bar{b}}(\bkR^n)\hookrightarrow Z_{p,q,n,\bar{b}}(\bkR^n)=:Z(\bkR^n),
\end{equation}
where
$$
Z(\bkR^n):=\left\{f\in{\mcal M}(\bkR^n):\|f\|_{Z}<+\infty\right\}
$$
and
\begin{equation}\label{eq:20}
\|f\|_{Z}=\|f\|_{Z_{p,q,n,\bar{b}}}:=\|t^{-1/q}\bar{b}(t^{1/n})\,\|f^*(y)\|_{p;(0,t)}\|_{q;(0,+\infty)}.
\end{equation}
\end{theorem}
\bpr
Putting
\begin{equation}\label{eq:21}
d(x):=\bar{b}(x^{p/n})\quad\text{for all}\quad x\in(0,+\infty),
\end{equation}
we see that $d\in SV(0,+\infty)$ and
\begin{align}
\|t^{-1/q} d(t)\|_{q;(0,1)}&\approx \|\tau^{-1/q}\bar{b}(\tau)\|_{q;(0,1)}=+\infty, \label{eq:22}\\
\|t^{-1/q} d(t)\|_{q;(1,+\infty)}&\approx \|\tau^{-1/q}\bar{b}(\tau)\|_{q;(1,+\infty)}<+\infty \label{eq:23}.
\end{align}
Using the trivial embedding $L_p(\bkR^n)\hookrightarrow L_p(\bkR^n)$ and \eqref{eq:17:embBintoLinfty}, we arrive at
$$
X(\bkR^n):=(L_p(\bkR^n), B_{p,1}^{n/p}(\bkR^n))_{0,q,d}\hookrightarrow (L_p(\bkR^n),L_{\infty}(\bkR^n))_{0,q,d}=:Z(\bkR^n).
$$
Choose $r\in\bkN$ such that $\theta:=\frac{n}{p}\frac{1}{r}<1$. Then $B_{p,1}^{n/p}(\bkR^n)=B_{p,1}^{\theta r}(\bkR^n)$. Together with 
\eqref{eq:16:interpLpW=B}, this yields
$$
X(\bkR^n)=(L_p(\bkR^n), (L_p(\bkR^n),W_p^{r}(\bkR^n))_{\theta,1})_{0,q,d},
$$
and, on applying \eqref{eq:3.20:GOT2002} of Theorem \ref{Theo3.5:GOT2002Lrrinwsvf}, with $X_0=L_p(\bkR^n)$ and $X_1=W_p^{r}(\bkR^n)$, we obtain
$$
X(\bkR^n)=(L_p(\bkR^n),W_p^{r}(\bkR^n))_{0,q,d\circ \rho}\quad\text{where}\quad \rho(t):=t^{\theta}\quad\text{for all}\quad t \in(0,+\infty).
$$
By \cite[Chap. V, Theorem 4.12, p. 339]{BS:IO},
$$
K(f,t;L_p(\bkR^n),W_p^{r}(\bkR^n)\approx \min(1,t)\|f\|_p+\omega_r(f,t^{1/r})_p
$$
for all $f\in L_p(\bkR^n)$ and all $t\in(0,+\infty)$. Thus, for all $f\in L_p(\bkR^n)$,
\begin{align}
\|f\|_X&\approx\|t^{-1/q}d(t^{\theta})\min(1,t)\|f\|_p\|_{q;(0,+\infty)} + \|t^{-1/q}d(t^{\theta})\,\omega_r(f,t^{1/r})_p\|_{q;(0,+\infty)}\label{eq:24}\\
&=:N_1+N_2.\nonumber
\end{align}
Moreover, $N_1\approx N_{11}+N_{12}$, where
$$
N_{11}:=\|f\|_p \cdot \|t^{1-1/q}d(t^{\theta})\|_{q;(0,1)}\approx \|f\|_p
$$
and, by \eqref{eq:23},
\begin{align}
N_{12}&:=\|f\|_p \cdot  \|t^{-1/q}d(t^{\theta})\|_{q;(1,+\infty)}\nonumber\\
&\approx \|f\|_p \cdot \|\tau^{-1/q}d(\tau)\|_{q;(1,+\infty)}\nonumber\\
&\approx \|f\|_p.\nonumber
\end{align}
Therefore,
\begin{equation}\label{eq:25}
N_1\approx \|f\|_p. 
\end{equation}
Furthermore, by \eqref{eq:21},
\begin{align}
N_{2}&:=\|t^{-1/q}d(t^{\theta})\,\omega_r(f,t^{1/r})_p\|_{q;(0,+\infty)}\nonumber\\
&\approx \|\tau^{-1/q}d(\tau^{\theta r})\,\omega_r(f,\tau)_p\|_{q;(0,+\infty)}\nonumber\\
&=\|\tau^{-1/q}d(\tau^{n/p})\,\omega_r(f,\tau)_p\|_{q;(0,+\infty)}\nonumber\\
&=\|\tau^{-1/q}\bar{b}(\tau)\,\omega_r(f,\tau)_p\|_{q;(0,+\infty)}.\nonumber
\end{align}
Together with \eqref{eq:24} and \eqref{eq:25}, this shows that the norm in the space $X(\bkR^n)$ is equivalent to
$$
\|\cdot\|_p +\|\tau^{-1/q}\bar{b}(\tau)\,\omega_r(\cdot,\tau)_p\|_{q;(0,+\infty)}.
$$
Hence, $X(\bkR^n)=B_{p,q}^{0,\bar{b}}(\bkR^n)$.
 
Now we are going to determine the space $Z(\bkR^n)$. Since, by \cite[Chap. V, Theorem~1.9, p. 300]{BS:IO}, $L_p(\bkR^n)=(L_1(\bkR^n), L_{\infty}(\bkR^n))_{\frac{1}{p'},p}$, we see that
$$
Z(\bkR^n)=((L_1(\bkR^n), L_{\infty}(\bkR^n))_{\frac{1}{p'},p}, L_{\infty}(\bkR^n))_{0,q,d}
$$
and on using \eqref{eq:3.21:GOT2002} of Theorem \ref{Theo3.5:GOT2002Lrrinwsvf}, with $X_0=L_1(\bkR^n)$ and $X_1=L_{\infty}(\bkR^n)$, we obtain
\begin{equation}\label{eq:26}
Z(\bkR^n)=\overline{X}_{\frac{1}{p'},q,d\circ\rho,p,1}^{\mcal L},\quad\text{where}\quad \rho(t)=t^{1/p},\quad t\in(0,+\infty).
\end{equation}
Hence,
$$
\|f\|_{Z}=\|t^{-1/q}d(t^{1/p})\,\|\tau^{-1/p'-1/p}\,K(f,\tau;L_1(\bkR^n),L_{\infty}(\bkR^n))\|_{p;(0,t)}\|_{q;(0,+\infty)}.
$$
Making use of the fact that
$
d(t^{1/p})=\bar{b}(t^{1/n})
$
for all $t\in (0,+\infty)$,
we arrive at$$
\|f\|_{Z}=\|t^{-1/q}\bar{b}(t^{1/n})\,\|\tau^{-1/p'-1/p}\,K(f,\tau;L_1(\bkR^n),L_{\infty}(\bkR^n))\|_{p;(0,t)}\|_{q;(0,+\infty)}.
$$
Finally, applying \eqref{eq:5.10:GOT2002} of Lemma \ref{lemma5.2:GOT2002Lrrinwsvf}, we obtain
$$
\|f\|_{Z}=\|t^{-1/q}\bar{b}(t^{1/n})\,\|f^{*}(y)\|_{p;(0,t)}\|_{q;(0,+\infty)}
$$
and the proof is complete.
\epr

\

\begin{proofspof}\ref{theo:mainresult}.
By Lemma \ref{Lemma10:aux},
\begin{equation}\label{embedding:lemmaaux:eq27}
B_{p,r}^{0,b}(\bkR^n)\hookrightarrow B_{p,q}^{0,\bar{b}}(\bkR^n),
\end{equation}
where
$$
\bar{b}(t):=(b_r(t))^{1-r/q}(b(t))^{r/q}\quad\text{for all}\quad t\in(0,+\infty).
$$
Using a change of variables and \eqref{cond:maintheo:1}, we get
\begin{align}
\|t^{-1/q}\bar{b}(t)\|_{q;(0,1)}&=\left(\int_0^1 t^{-1}\left(\int_t^{+\infty}\tau^{-1}b^r(\tau)\,d\tau\right)^{\frac{q}{r}-1}\,b^r(t)\,d t \right)^{1/q}\label{tag16*}\\
&=\left(\int_{(b_r(1))^r}^{+\infty}y^{\frac{q}{r}-1}\,dy\right)^{1/q}=+\infty\nonumber
\end{align}
and
\begin{equation}\label{tag16**}
\|t^{-1/q}\bar{b}(t)\|_{q;(1,+\infty)}=\left(\int_0^{(b_r(1))^r}y^{\frac{q}{r}-1}\,dy\right)^{1/q}\approx b_r(1)<+\infty.
\end{equation}
Therefore, by Theorem \ref{theorem11},
\begin{equation}\label{eq:28}
B_{p,q}^{0,\bar{b}}(\bkR^n)\hookrightarrow Z_{p,q,n,\bar{b}}(\bkR^n),
\end{equation}
where
$$
\|f\|_{Z}=\|t^{-1/q}\bar{b}(t^{1/n})\,\|f^*(y)\|_{p;(0,t)}\|_{q;(0,+\infty)}.
$$
Combining now \eqref{embedding:lemmaaux:eq27} and \eqref{eq:28}, we obtain the result.
\hfill $\square$ 
\end{proofspof}

\begin{remark}\label{remark12}
{\rm Note that if $n\in\bkN$, $1\leq r<+\infty$ and $a\in SV(0,+\infty)$, then, for all $t\in(0,+\infty)$,
\begin{align*}
a_r(t^{1/n})&:=\|\tau^{-1/r}a(\tau)\|_{r;(t^{1/n},+\infty)}=\left(\int_{t^{1/n}}^{+\infty} \tau^{-1}a^r(\tau)\,d\tau \right)^{1/r}\\
&\approx \left(\int_{t}^{+\infty} y^{-1}a^r(y^{1/n})\,d y\right)^{1/r}=:a_{r,n}(t),
\end{align*}
i.e.,
\begin{equation}\label{eq:29}
a_r(t^{1/n})\approx a_{r,n}(t)\quad\text{for all}\quad t\in(0,+\infty).
\end{equation}
}
\end{remark}

\

\begin{proofoff}\ref{cor:mainresult}.
By the sufficient part of Theorem \ref{theo:mainresult}, embedding \eqref{embedding:maintheo:1} holds, with $\bar{b}$ given by \eqref{deffunct:maintheo:1}.

If $d$ is defined by \eqref{eq:21}, then, by \eqref{eq:26},
\begin{equation}\label{eq:30aux}
Z(\bkR^n)=\overline{X}_{\frac{1}{p'},q,d\circ\rho,p,1}^{\mcal L},\quad\text{where}\quad \rho(t)=t^{1/p},\quad t\in (0,+\infty).
\end{equation}
Thus, by Theorem \ref{Theorem5:circ},
\begin{equation}\label{eq:30}
Z(\bkR^n)\hookrightarrow \overline{X}_{\frac{1}{p'},q,\widetilde{d}}:=(L_1(\bkR^n),L_{\infty}(\bkR^n))_{\frac{1}{p'},q,\widetilde{d}}\,,
\end{equation}
where
\begin{align}
\widetilde{d}(t)&:= (d\circ\rho)(t)\left(\frac{\int_t^{+\infty}y^{-1}\left[(d\circ\rho)(y)\right]^q\,dy}{\left[(d\circ\rho)(t)\right]^q}\right)^{\frac{1}{\max\{p,q\}}}\label{eq:31}\\
&=d(t^{1/p})\left(\frac{\|y^{-1/q} d(y^{1/p})\|_{q; (t,+\infty)}}{d(t^{1/p})}\right)^{\frac{q}{\max\{p,q\}}},\quad t\in(0,+\infty).\nonumber
\end{align}
Since $d(t^{1/p})=\bar{b}(t^{1/n})$ for any $t\in(0,+\infty)$, \eqref{eq:31} can be rewritten as
\begin{align}
\widetilde{d}(t)&:= \bar{b}(t^{1/n})\left(\frac{\|y^{-1/q} \bar{b}(y^{1/n})\|_{q; (t,+\infty)}}{\bar{b}(t^{1/n})}\right)^{\frac{q}{\max\{p,q\}}}\label{eq:32}\\
&=\bar{b}(t^{1/n})\left(\frac{\bar{b}_{q,n}(t)}{\bar{b}(t^{1/n})}\right)^{\frac{q}{\max\{p,q\}}},\quad t\in(0,+\infty).\nonumber
\end{align}
First we calculate $\bar{b}_{q,n}$. Using \eqref{eq:29}, \eqref{deffunct:maintheo:1} and a change of variables, we obtain, for all $t\in(0,+\infty)$,
\begin{align*}
\bar{b}_{q,n}(t)&\approx \bar{b}_{q}(t^{1/n})=\left(\int_{t^{1/n}}^{+\infty}y^{-1}\bar{b}^q(y)\,dy\right)^{1/q}\\
&=\left(\int_{t^{1/n}}^{+\infty}y^{-1}\left(\int_y^{+\infty}\tau^{-1}b^r(\tau)\,d\tau\right)^{\frac{q}{r}-1}b^r(y)\,dy\right)^{1/q}\\
&\approx\left(\int_0^{(b_r(t^{1/n}))^r}x^{\frac{q}{r}-1}\,dx\right)^{1/q}\\
&\approx\left(\int_{t^{1/n}}^{+\infty}\tau^{-1}b^r(\tau)\,d\tau\right)^{1/r}=b_r(t^{1/n}),
\end{align*}
i.e., 
\begin{equation}\label{eq:33:aux}
\bar{b}_{q,n}(t)\approx \bar{b}_q(t^{1/n}) \approx b_r(t^{1/n})\quad\text{for all}\quad t\in(0,+\infty).
\end{equation}
Thus,
$$\widetilde{d}(t)\approx\bar{b}(t^{1/n})\left(\frac{b_r(t^{1/n})}{\bar{b}(t^{1/n})}\right)^{\frac{q}{\max\{p,q\}}}\quad\text{for all}\quad t\in(0,+\infty),
$$
and, on using \eqref{deffunct:maintheo:1} and \eqref{deffunct:maintcor:1}, we arrive at
\begin{align}\label{eq:33NEW}
\widetilde{d}(t)&\approx (b_r(t^{1/n}))^{1-\frac{r}{q}}(b(t^{1/n}))^{\frac{r}{q}}\left(\frac{b_r(t^{1/n})}{(b_r(t^{1/n}))^{1-\frac{r}{q}}(b(t^{1/n}))^{\frac{r}{q}}}\right)^{\frac{q}{\max\{p,q\}}}\\
&=(b_r(t^{1/n}))^{1-\frac{r}{q}+\frac{r}{\max\{p,q\}}}(b(t^{1/n}))^{\frac{r}{q}-\frac{r}{\max\{p,q\}}}\nonumber\\
&=\widetilde{b}(t)\quad\text{for all}\quad  t\in(0,+\infty).\nonumber
\end{align}
Therefore, by \eqref{embedding:maintheo:1} and \eqref{eq:30}, 
\begin{equation}\label{eq:33NEWaux}
B_{p,r}^{0,b}(\bkR^n)\hookrightarrow (L_1(\bkR^n),L_{\infty}(\bkR^n))_{\frac{1}{p'},q,\widetilde{d}}\,.
\end{equation}
Since one can easily show that
\begin{equation}\label{eq:33aux}
(L_1(\bkR^n),L_{\infty}(\bkR^n))_{\frac{1}{p'},q,\widetilde{d}}=L_{p,q;\widetilde{d}}\;(\bkR^n),
\end{equation}
together with \eqref{eq:33NEW} and \eqref{eq:33NEWaux}, this gives \eqref{embedding:maincor:1}.
\hfill $\square$ 
\end{proofoff}

\

\begin{proofnpof}\ref{theo:mainresult}. 
Assume now that \eqref{embedding:maintheo:1} holds. Then, as in the proof of Corollary \ref{cor:mainresult}, one can show that
$$
Z_{p,q,n,\bar{b}}(\bkR^n)\hookrightarrow L_{p,q;\widetilde{b}}(\bkR^n),
$$
with $\widetilde{b}$ defined by \eqref{deffunct:maintcor:1}. Consequently,
$$
B_{p,r}^{0,b}(\bkR^n)\hookrightarrow L_{p,q;\widetilde{b}}^{loc}(\bkR^n).
$$
Then, together with \eqref{eq:29},  Theorem \ref{main*:CGO2} gives $q\geq r$.
\hfill $\square$ 
\end{proofnpof}

%
%
%

\section{Proofs of Sharpness and Local Optimality}


\begin{lemma}\label{lemma14}
Let $n\in\bkN$, $1<p<+\infty$, $1\leq r\leq q<+\infty$, and let $b\in SV(0,+\infty)$ satisfy \eqref{cond:maintheo:1}. If $\,\bar{b}_{(r,r)}, \bar{b}_{(r,q)}\in SV(0,+\infty)$ are given by  \eqref{deffunct:maintheo:1}, then 
\begin{equation}\label{eq:35}
Z_{p,r,n,b}(\bkR^n)= Z_{p,r,n,\bar{b}_{(r,r)}}(\bkR^n)\hookrightarrow Z_{p,q,n,\bar{b}_{(r,q)}}(\bkR^n).
\end{equation}
\end{lemma}
\bpr
Embedding \eqref{eq:35} means that, for all $f\in Z_{p,r,n,\bar{b}_{(r,r)}}(\bkR^n)$, \footnotemark$^{)}$\footnotetext{$^)$ Note that $\bar{b}_{(r,r)}=b$ and so $Z_{p,r,n,\bar{b}_{(r,r)}}(\bkR^n)= Z_{p,r,n,b}(\bkR^n)$.}
\begin{align*}
\big\|t^{-1/q}\bar{b}_{(r,q)}(t^{1/n})\,\|f^*\|_{p;(0,t)}\big\|_{q;(0,+\infty)}&\precsim
\big\|t^{-1/r}\bar{b}_{(r,r)}(t^{1/n})\,\|f^*\|_{p;(0,t)}\big\|_{r;(0,+\infty)}\\
&=\big\|t^{-1/r}b(t^{1/n})\,\|f^*\|_{p;(0,t)}\big\|_{r;(0,+\infty)},
\end{align*}
which can be rewritten as
\begin{equation}\label{eq:36}
\left(\int_0^{+\infty}w(t)g^q(t)\,dt\right)^{1/q}\precsim \left(\int_0^{+\infty}v(t)g^r(t)\,dt\right)^{1/r},
\end{equation}
where, for all $t\in(0,+\infty)$,
$$
g(t):=\|f^*\|_{p;(0,t)},\quad w(t):=t^{-1}{\bar{b}_{(r,q)}}^q(t^{1/n}),\quad\text{and}\quad v(t):=t^{-1}{\bar{b}_{(r,r)}}^r(t^{1/n})=t^{-1}b^r(t^{1/n}).
$$
Since $g\in\mathcal{M}^+(0,+\infty;\uparrow)$, on putting
$$
W(t):=\int_t^{+\infty}w(\tau)\,d\tau\quad\text{and}\quad V(t):=\int_t^{+\infty}v(\tau)\,d\tau\quad\text{for all}\quad t\in(0,+\infty),
$$
we obtain by \cite[Proposition 2.1]{H-Step} that inequality \eqref{eq:36} holds for all $g\in\mathcal{M}^+(0,+\infty;\uparrow)$ if and only if
\begin{equation}\label{eq:37}
\sup_{t\in(0,+\infty)}W^{1/q}(t)V^{-1/r}(t)<+\infty.
\end{equation}
As $V^{1/r}(t)\approx \bar{b}_{r}(t^{1/n})$ and $W^{1/q}(t)\approx \bar{b}_{r}(t^{1/n})$ for all $t\in(0,+\infty)$, condition \eqref{eq:37} is satisfied. Consequently, \eqref{eq:35} holds.
\epr

\bprof \ref{theo:relatedtolemma14}.
The result is an immediate consequence of Theorem \ref{theo:mainresult} and Lemma \ref{lemma14}. 
\eprof

\bprof \ref{sharpnesswithrespecttofunction}. Putting $\bar{B}=\kappa\bar{b}$ (with $\bar{b}$ given by \eqref{deffunct:maintheo:1}), then one can show (as in the proof of Corollary \ref{cor:mainresult})  that
\begin{equation}\label{em:aux:ZL}
Z_{p,q,n,\bar{B}}(\bkR^n)\hookrightarrow L_{p,q;\widetilde{B}}(\bkR^n),
\end{equation}
where (cf. \eqref{eq:32}  and \eqref{eq:29}), for all $t\in(0,+\infty)$,
\begin{align}
\widetilde{B}(t)
&=\bar{B}(t^{1/n})\left(\frac{\bar{B}_{q,n}(t)}{\bar{B}(t^{1/n})}\right)^{\frac{q}{\max\{p,q\}}}\label{tag19*}\\
&\approx \bar{B}(t^{1/n})\left(\frac{\bar{B}_{q}(t^{1/n})}{\bar{B}(t^{1/n})}\right)^{\frac{q}{\max\{p,q\}}}\nonumber\\
&= \kappa(t^{1/n})\bar{b}(t^{1/n})\left(\frac{\|\tau^{-1/q}\kappa(\tau)\bar{b}(\tau)\|_{q;({t^{1/n},+\infty})}}{\kappa(t^{1/n})\bar{b}(t^{1/n})}\right)^{\frac{q}{\max\{p,q\}}}\nonumber.
\end{align}
Hence, by \eqref{eq:41} and \eqref{em:aux:ZL},
\begin{equation*}\label{eq:41:aux1}
\big\|t^{1/p-1/q}\widetilde{B}(t)f^*\big\|_{q;(0,+\infty)}\precsim \|f\|_{B_{p,r}^{0,b}}\quad\text{for all}\quad f\in B_{p,r}^{0,b}(\bkR^n).\end{equation*}
In particular, this gives
\begin{equation}\label{eq:41:aux2}
\big\|t^{1/p-1/q}\widetilde{B}(t)f^*(t)\big\|_{q;(0,1)}\precsim \|f\|_{B_{p,r}^{0,b}}\quad\text{for all}\quad f\in B_{p,r}^{0,b}(\bkR^n).
\end{equation}
Now we intend to use Theorem \ref{T1:CGO2} to verify the validity of \eqref{eq:41:aux2}. To this end, we are going to 
estimate the quantity
$$\big\|s^{-1/p-1/\rho}{\mcal W}_q(s)\big\|_{
\rho;(t,1)}\quad\text{for all}\quad t\in(0,1),
$$
where $${\mcal W}_q(t):=\big\|s^{1/p-1/q}\widetilde{B}(s)\big\|_{q;(0,t)}, \quad t\in(0,1],$$ 
and $\rho=+\infty$ if $p\leq q$, and $\rho$ is defined by $\frac{1}{\rho}=\frac{1}{q}-\frac{1}{p}$ if $q<p$.
Since $\widetilde{B}\in SV(0,+\infty)$,
$${\mcal W}_q(t)\approx t^{1/p}\widetilde{B}(t)\quad\text{for all}\quad t\in(0,1),$$
which implies that
 $$\big\|s^{-1/p-1/\rho}{\mcal W}_q(s)\big\|_{\rho;(t,1)}\approx\big\|s^{-1/\rho}\widetilde{B}(s)\big\|_{\rho;(t,1)}\quad\text{for all}\quad t\in(0,1).$$
If $p\leq q$, then, by \eqref{tag19*},
 $$\big\|s^{-1/\rho}\widetilde{B}(s)\big\|_{\rho;(t,1)}=\|\bar{B}_{q,n}(s)\big\|_{\infty;(t,1)}\approx \bar{B}_{q}(t^{1/n})\quad\text{for all}\quad t\in(0,1).$$
 If $q< p$, then, by \eqref{tag19*} and a simple calculation,
 \begin{align*}
 \big\|s^{-1/\rho}\widetilde{B}(s)\big\|_{\rho;(t,1)}&\approx\big\|s^{-1/\rho}(\bar{B}(s^{1/n}))^{1-q/p}(\bar{B}_{q}(s^{1/n}))^{q/p}\big\|_{\rho;(t,1)}\\
& \approx\big\|s^{-1/\rho}(\bar{B}(s^{1/n}))^{q/\rho}(\bar{B}_{q}(s^{1/n}))^{1-q/\rho}\big\|_{\rho;(t,1)}\\
 &\approx((\bar{B}_{q}(t^{1/n}))^{\rho}-(\bar{B}_{q}(1))^{\rho})^{1/\rho}
 \quad\text{for all}\quad t\in(0,1).
 \end{align*}
Put $A=0$ if  $\rho=+\infty$ and $A=(\bar{B}_{q}(1))$ if $\rho<+\infty$.
First, using L'Hopital's Rule, one can see that
 $$
 \lim_{t\rightarrow 0_+}\frac{(\bar{B}_{q}(t^{1/n}))^{\rho}-A^{\rho})^{1/\rho}}{\bar{b}_{q}(t^{1/n})}=\kappa(0_+).
 $$
Together with \eqref{eq:33:aux} and \eqref{eq:29}, this implies that 
\eqref{T1.1} cannot hold unless $\kappa(0_+)<+\infty$, that is, unless $\kappa$ is bounded.
\eprof

\section{Relations among the target spaces in the main results}

\bprof \ref{theorem24}.
First, by \eqref{eq:30}, \eqref{eq:33NEW} and \eqref{eq:33aux}, 
$$
Z_{p,q,n,\bar{b}}(\bkR^n)\hookrightarrow  L_{p,q;\widetilde{b}}(\bkR^n).
$$
Moreover, by \eqref{eq:30aux}, 
$$
Z_{p,q,n,\bar{b}}(\bkR^n)=\overline{X}_{\frac{1}{p'},q,d\circ\rho,p,1}^{\mcal L}.
$$
Hence, it follows from \eqref{eq:33NEW}, \eqref{eq:33aux} and \eqref{eq:56} that $Z_{p,q,n,\bar{b}}(\bkR^n)=L_{p,q;\widetilde{b}}(\bkR^n)$ if $q=p$. 

Now we intend to prove that the opposite embedding to \eqref{eq:61aux1} does not hold if $q\neq p$. That is, we are going to prove
\eqref{eq:62}, which is
$$
 L_{p,q;\widetilde{b}}(\bkR^n)\hookrightarrow{\!\!\!\!\!\!\!/}\,\;\;Z_{p,q,n,\bar{b}}(\bkR^n)\quad\text{if}\quad q\neq p.
$$

%
On the contrary, suppose that 
\begin{equation}\label{eq:63}
 L_{p,q;\widetilde{b}}(\bkR^n)\hookrightarrow Z_{p,q,n,\bar{b}}(\bkR^n)
\end{equation}
if $p\neq q$, $1<p<+\infty$, $1\leq r\leq q<+\infty$ and $b\in SV(0,+\infty)$ safisfies \eqref{cond:maintheo:1}.\footnotemark$^{)}$\footnotetext{$^{)}$ Recall that $r$ is involved in the definitions of $\bar{b}$ and $\widetilde{b}$.}
Then, according to the proof of Theorem \ref{theo:mainresult} (cf. \eqref{tag16*} and  \eqref{tag16**}), $\bar{b}$ satisfies
\begin{equation}\label{eq:63half}
\|t^{-1/q}\bar{b}(t)\|_{q;(0,1)}=+\infty\quad\text{and}\quad \|t^{-1/q}\bar{b}(t)\|_{q;(1,+\infty)}<+\infty.\mbox{\footnotemark$^{)}$}
\end{equation}
\footnotetext{$^{)}$ Note that \eqref{eq:63half} implies that $0\neq\bar{b}(t)\neq +\infty$ for all $t\in(0,+\infty)$.
}

Embeding \eqref{eq:63} means that 
$$
\|f\|_{Z_{p,q,n,\bar{b}}}\precsim \|f\|_{L_{p,q;\widetilde{b}}}\quad\text{for all}\quad f \in L_{p,q;\widetilde{b}}(\bkR^n),
$$
which is equivalent to 
$$
\left\|t^{-\frac{1}{q}}\bar{b}(t^{\frac{1}{n}}\!)\!\left(\int_0^t\!\!(f^*(y))^p\,dy\!\right)^{\frac{1}{p}}\right\|_{q;(0,+\infty)}
\!\!\precsim 
\|t^{\frac{1}{p}-\frac{1}{q}}\widetilde{b}(t)f^*(t)\|_{q;(0,+\infty)}
\text{ \ for all \ } f \in L_{p,q;\widetilde{b}}(\bkR^n).
$$
Hence, putting $g:=(f^*)^p$ and using \cite[Chap. III, Corollary 7.8, p. 86]{BS:IO}, we obtain
\begin{equation}\label{eq:64}
\left\|t^{-\frac{1}{q}}\bar{b}(t^{\frac{1}{n}}\!)\!\left(\int_0^t\!\!g(y)\,dy\!\right)^{\frac{1}{p}}\right\|_{q;(0,+\infty)}
\!\!\!\!\!\!\!\precsim 
\|t^{\frac{1}{p}-\frac{1}{q}}\widetilde{b}(t)g^{\frac{1}{p}}(t)\|_{q;(0,+\infty)}
\text{ \ for all \ }g \!\in\!{\mcal M}_0^+(0,\!+\infty;\downarrow).
\end{equation}
Set $Q:=\frac{q}{p}$. Then \eqref{eq:64} is equivalent to
\begin{align}\label{eq:65}
\left(\int_0^{+\infty}\!\!\!t^{-1}\bar{b}^q(t^{\frac{1}{n}}\!)\right.&\left.\!\left(\int_0^t\!\!g(y)dy\!\right)^{Q}\!\!dt\right)^{1/Q}\\
&\precsim 
\left(\int_0^{+\infty}\!\!\! t^{Q-1}\widetilde{b}^q(t)g^{Q}(t)dt\right)^{1/Q}\quad\text{for all}\quad
g \!\in\!{\mcal M}_0^+(0,\!+\infty;\downarrow).\nonumber
\end{align}
\begin{trivlist}
\item[\hspace*{0.5cm}{\rm I.}] Assume first that $Q=\frac{q}{p}\in (0,1)$, i.e., $q<p$. Inequality \eqref{eq:65} can be rewritten as
\begin{align}\label{eq:66}
\left(\int_0^{+\infty}\right.&\left.\left(\int_0^x\!\!g(y)dy\!\right)^{Q}\!w(x)dx\right)^{1/Q}\\
&\precsim 
\left(\int_0^{+\infty} \!\!\!g^{Q}(x) v(x)\,dx\right)^{1/Q}\quad\text{for all}\quad g \!\in\!{\mcal M}_0^+(0,\!+\infty;\downarrow),\nonumber
\end{align}
where, for all $x\in (0,+\infty)$,
$$
w(x):=x^{-1}\,\bar{b}^q(x^{1/n})\quad\text{and}\quad
v(x):=x^{Q-1}\,\widetilde{b}^q(x).
$$
If we put $P:=Q$ and
$$
\varphi(x,y):=\chi_{(0,x)}(y)\quad\text{for all}\quad x,y\in(0,+\infty),
$$
then \eqref{eq:66} coincides with inequality \eqref{eq7:circ:STAR}. This fact and Theorem \ref{theorem3=Theorem2.1:Lai:1993wt:STAR} then imply that
\begin{equation}\label{eq:67}
 \left(\int_0^{+\infty}\left(\int_0^{z} \varphi(x,y) \,dy\right)^Q\,w(x)\,dx\right)^{1/Q}\precsim \left(\int_0^{z}v(x)\,dx\right)^{1/Q}
\text{ \ for all \ } z\in(0,+\infty).
\end{equation}
However, we are going to show that this is not the case.
First,
\begin{equation}\label{eq:68}
 {\rm LHS \eqref{eq:67}} \approx z \left(\bar{b}_q(z^{1/n})\right)^p
 \quad\text{ for all}\quad z\in(0,+\infty).
\end{equation}
Second,
\begin{equation}\label{eq:69}
 {\rm RHS \eqref{eq:67}} \approx z \left(\widetilde{b}(z)\right)^p
 \quad\text{ for all}\quad z\in(0,+\infty).
\end{equation}
Since \eqref{deffunct:maintheo:1} implies that (cf. \eqref{eq:32}, \eqref{eq:33:aux} and \eqref{eq:33NEW})
\begin{equation}\label{eq:70}
\widetilde{b}(z)\approx\bar{b}(z^{1/n})\left(\frac{\bar{b}_q(z^{1/n})}{\bar{b}(z^{1/n})}\right)^{\frac{q}{\max\{p,q\}}}\quad\text{for all}\quad z\in(0,+\infty),
\end{equation}
we see that
\begin{equation}\label{eq:72}
(\widetilde{b}(z))^p\approx \left(\bar{b}(z^{1/n})\right)^{p-q} \left(\bar{b}_q(z^{1/n})\right)^q  \quad\text{ for all}\quad z\in(0,+\infty).
\end{equation}
Therefore, by \eqref{eq:68}, \eqref{eq:69} and \eqref{eq:72}, inequality \eqref{eq:67} is equivalent to
$$
\left(\bar{b}_q(z^{1/n})\right)^p\precsim \left(\bar{b}(z^{1/n})\right)^{p-q} \left(\bar{b}_q(z^{1/n})\right)^q  \quad\text{ for all}\quad z\in(0,+\infty),
$$
which can be rewritten as
$$
\left(\int_{z^{1/n}}^{+\infty}\tau^{-1}\,\bar{b}^q(\tau)\,d\tau\right)^{1/q}\precsim \bar{b}(z^{1/n})  \quad\text{ for all}\quad z\in(0,+\infty).
$$
However, this is a contradiction since, by \eqref{prop:SS_Limite_Frac_0},
$$
\limsup_{x\rightarrow 0_+}\dfrac{\|\tau^{-1/q} \bar{b}(\tau)\|_{q;(x,+\infty)}}{\bar{b}(x)}=+\infty.
$$

\item[\hspace*{0.5cm}{\rm II.}] Assume now that $Q=\frac{q}{p}\in (1,+\infty)$, i.e., $p<q$. Inequality \eqref{eq:65} can be rewritten as
\eqref{eq:58} with $P=Q$ and
$$
v(x):=x^{Q-1}\,(\widetilde{b}(x))^q\quad\text{and}\quad
w(x):=x^{Q-1}\,\bar{b}^q(x^{1/n}) \quad\text{ for all}\quad x\in (0,+\infty).
$$
This fact and Theorem \ref{theorem23=theorem2Sawyer} then imply that both conditions \eqref{eq:59} and \eqref{eq:60} are satisfied. However, we are going to show that this is not the case. 

Let us check condition \eqref{eq:60}. First,
\begin{align}
\left(\int_z^{+\infty}x^{-Q}\,w(x)\,dx\right)^{1/Q}&=\left(\int_z^{+\infty}x^{-1}\,\bar{b}^q(x^{1/n})\,dx\right)^{1/Q}\label{eq:73}\\
&=\left(\bar{b}_{q,n}(z)\right)^p\approx \left(\bar{b}_q(z^{1/n})\right)^p\quad\text{ for all}\quad z\in(0,+\infty).\nonumber
\end{align}

\noindent Second,
\begin{align}
\left(\!\int_0^{z}\!\!\left(\frac{x}{\int_0^{x}\!\!v(t)\,dt}\!\right)^{Q'}\!\!\!v(x)\,dx\!\right)^{1/Q'}\!\!&=
\left(\!\int_0^z\!\! x^{Q'}\!\!\left(\int_0^x\!\!t^{Q-1}(\widetilde{b}(t))^q\,\!dt\!\right)^{-Q'}\!\!\!x^{Q-1}\,(\widetilde{b}(x))^q\,dx\!\right)^{1/Q'}\label{eq:74}\\
&\approx \left(\int_0^z\left(x^{Q-1}(\widetilde{b}(x))^q\right)^{-Q'}x^{Q-1}\,(\widetilde{b}(x))^q\,dx\right)^{1/Q'}\nonumber\\
&=\left(\int_0^z x^{-1}\left((\widetilde{b}(x))^q\right)^{1-Q'}\,dx\right)^{1/Q'}\nonumber\\
&=\left(\int_0^z x^{-1}\left((\widetilde{b}(x))^{-\frac{q}{Q}}\right)^{Q'}\,dx\right)^{1/Q'}\nonumber\\
&=\left(\int_0^z x^{-1}\left((\widetilde{b}(x))^{-p}\right)^{Q'}\,dx\right)^{1/Q'}\quad\text{ for all}\quad z>0.\nonumber
\end{align}
Since now $q>p$,  \eqref{eq:70} 
implies that
$$
(\widetilde{b}(z))^{-p}\approx \left(\bar{b}(z^{1/n})\dfrac{\bar{b}_{q}(z^{1/n})}{\bar{b}(z^{1/n})}\right)^{-p} = \left(\bar{b}_{q}(z^{1/n})\right)^{-p}\quad\text{ for all}\quad z\in(0,+\infty).
$$
Thus, defining the function $a\in SV(0,+\infty)$ by
$$
a(x):=\left(\bar{b}_{q}(x^{1/n})\right)^{-p},\quad x\in(0,+\infty),
$$
we obtain from \eqref{eq:73} and \eqref{eq:74} that
$$
 {\rm LHS \eqref{eq:60}} \approx \dfrac{\left(\int_0^z x^{-1}a^{Q'}(x)\,dx\right)^{1/Q'}}{a(z)}
 \quad\text{ for all}\quad z\in(0,+\infty).
$$
Since, by \eqref {prop:SS_Limite_Frac_Infty},
$$
\limsup_{z\rightarrow +\infty} \dfrac{\left(\int_0^zx^{-1}a^{Q'}(x)\,dx\right)^{1/Q'}}{a(z)}=+\infty,
$$
condition  \eqref{eq:60} does not hold, which is a contradiction.

\end{trivlist}

Consequently, in both cases I and II, \eqref{eq:62} holds.
\eprof

\bibliographystyle{alpha}

\begin{thebibliography}{CDSFM15}

\bibitem[BC18]{BeCo}
B.~F. Besoy and F.~Cobos.
\newblock Duality for logarithmic interpolation spaces when {$0<q<1$} and
  applications.
\newblock {\em J. Math. Anal. Appl.}, 466(1):373--399, 2018.

\bibitem[BGT87]{BGT87:RV}
N.~H. Bingham, C.~M. Goldie, and J.~L. Teugels.
\newblock {\em Regular Variation}.
\newblock Cambridge University Press, Cambridge, 1987.

\bibitem[BR80]{BR80:OLZS}
C.~Bennett and K.~Rudnick.
\newblock On {L}orentz-{Z}ygmund spaces.
\newblock {\em Dissertationes Math. (Rozprawy Mat.)}, 175:1--72, 1980.

\bibitem[BS88]{BS:IO}
C.~Bennett and R.~Sharpley.
\newblock {\em Interpolation of Operators}, volume 129 of {\em Pure and Applied
  Mathematics}.
\newblock Academic Press, New York, 1988.

\bibitem[CD14]{CD1}
F.~Cobos and {\'O}.~Dom\'{i}nguez.
\newblock Embeddings of {B}esov spaces of logarithmic smoothness.
\newblock {\em Studia Math.}, 223(3):193--204, 2014.

\bibitem[CD15a]{COBOS201543}
F.~Cobos and {\'O}.~Dom\'{\i}nguez.
\newblock Approximation spaces, limiting interpolation and {B}esov spaces.
\newblock {\em J. Approx. Theory}, 189:43 -- 66, 2015.

\bibitem[CD15b]{CD2}
F.~Cobos and {\'O}.~Dom\'{i}nguez.
\newblock On {B}esov spaces of logarithmic smoothness and {L}ipschitz spaces.
\newblock {\em J. Math. Anal. Appl.}, 425(1):71--84, 2015.

\bibitem[CD16]{CD3}
F.~Cobos and {\'O}.~Dom\'{i}nguez.
\newblock On the relationship between two kinds of {B}esov spaces with
  smoothness near zero and some other applications of limiting interpolation.
\newblock {\em J. Fourier Anal. Appl.}, 22(5):1174--1191, 2016.

\bibitem[CDK18]{CDK2018}
F.~Cobos, {\'O}.~Dom\'{i}nguez, and T.~K\"{u}hn.
\newblock Approximation and entropy numbers of embeddings between approximation
  spaces.
\newblock {\em Constr. Approx.}, 47(3):453--486, 2018.

\bibitem[CDSFM15]{CDSFM2015}
F.~Colombini, D.~Del~Santo, F.~Fanelli, and G.~M\'{e}tivier.
\newblock The well-posedness issue in {S}obolev spaces for hyperbolic systems
  with {Z}ygmund-type coefficients.
\newblock {\em Comm. Partial Differential Equations}, 40(11):2082--2121, 2015.

\bibitem[CDT16]{CDT}
F.~Cobos, {\'O}.~Dom\'{i}nguez, and H.~Triebel.
\newblock Characterizations of logarithmic {B}esov spaces in terms of
  differences, {F}ourier-analytical decompositions, wavelets and semi-groups.
\newblock {\em J. Funct. Anal.}, 270(12):4386--4425, 2016.

\bibitem[CF05]{CF2005}
C.~Capone and A.~Fiorenza.
\newblock On small {L}ebesgue spaces.
\newblock {\em J. Funct. Spaces Appl.}, 3(1):73--89, 2005.

\bibitem[CF06]{CaFa04}
A.~M. Caetano and W.~Farkas.
\newblock Local growth envelopes of {B}esov spaces of generalized smoothness.
\newblock {\em Z. Anal. Anwend.}, 25(3):265--298, 2006.

\bibitem[CGO08]{CGO1}
A.~M. Caetano, A.~Gogatishvili, and B.~Opic.
\newblock Sharp embeddings of {B}esov spaces involving only logarithmic
  smoothness.
\newblock {\em J. Approx. Theory}, 152(2):188--214, 2008.

\bibitem[CGO11a]{CGO3}
A.~Caetano, A.~Gogatishvili, and B.~Opic.
\newblock Compact embeddings of {B}esov spaces involving only slowly varying
  smoothness.
\newblock {\em Czechoslovak Math. J.}, 61(136)(4):923--940, 2011.

\bibitem[CGO11b]{CGO2}
A.~M. Caetano, A.~Gogatishvili, and B.~Opic.
\newblock Embeddings and the growth envelope of {B}esov spaces involving only
  slowly varying smoothness.
\newblock {\em J. Approx. Theory}, 163(10):1373--1399, 2011.

\bibitem[CH05]{CaHa04}
A.~M. Caetano and D.~D. Haroske.
\newblock Continuity envelopes of spaces of generalised smoothness: a limiting
  case; embeddings and approximation numbers.
\newblock {\em J. Funct. Spaces Appl.}, 3(1):33--71, 2005.

\bibitem[CL13]{CL2013}
A.~M. Caetano and H.-G. Leopold.
\newblock On generalized {B}esov and {T}riebel-{L}izorkin spaces of regular
  distributions.
\newblock {\em J. Funct. Anal.}, 264(12):2676--2703, 2013.

\bibitem[CM04a]{CM04b}
A.~M. Caetano and S.~D. Moura.
\newblock Local growth envelopes of spaces of generalized smoothness: the
  critical case.
\newblock {\em Math. Inequal. Appl.}, 7(4):573--606, 2004.

\bibitem[CM04b]{CM04a}
A.~M. Caetano and S.~D. Moura.
\newblock Local growth envelopes of spaces of generalized smoothness: the
  subcritical case.
\newblock {\em Math. Nachr.}, 273:43--57, 2004.

\bibitem[Dom16]{D2016}
{\'O}.~Dom\'{i}nguez.
\newblock Tractable embeddings of {B}esov spaces into small {L}ebesgue spaces.
\newblock {\em Math. Nachr.}, 289(14-15):1739--1759, 2016.

\bibitem[Dom17a]{DomO17}
{\'O}.~Dom\'{i}nguez.
\newblock Sharp embeddings of {B}esov spaces with logarithmic smoothness in
  sub-critical cases.
\newblock {\em Anal. Math.}, 43(2):219--240, 2017.

\bibitem[Dom17b]{D2017}
{\'O}.~Dom\'{i}nguez.
\newblock Ul'yanov-type inequalities and embeddings between {B}esov spaces: the
  case of parameters with limit values.
\newblock {\em Math. Inequal. Appl.}, 20(3):755--772, 2017.

\bibitem[DRS79]{DeVRS}
R.~A. DeVore, S.~D. Riemenschneider, and R.~C. Sharpley.
\newblock Weak interpolation in {B}anach spaces.
\newblock {\em J. Funct. Anal.}, 33(1):58--94, 1979.

\bibitem[EE04]{EEv04:HOFSE}
D.~E. Edmunds and W.~D. Evans.
\newblock {\em Hardy {O}perators, {F}unction {S}paces and {E}mbeddings}.
\newblock Springer-Verlag, Berlin, Heidelberg, 2004.

\bibitem[EGO97]{EGO4}
D.~E. Edmunds, P.~Gurka, and B.~Opic.
\newblock On embeddings of logarithmic \protect{B}essel potential spaces.
\newblock {\em J. Funct. Anal.}, 146(1):116--150, 1997.

\bibitem[EKP00]{EKP:OISRIQ}
D.~E. Edmunds, R.~Kerman, and L.~Pick.
\newblock Optimal {S}obolev {I}mbeddings {I}nvolving
  {R}earrangement-{I}nvariant {Q}uasinorms.
\newblock {\em J. Funct. Anal.}, 170:307--355, 2000.

\bibitem[EO00]{EO00:RILFR}
W.~D. Evans and B.~Opic.
\newblock Real interpolation with logarithmic functors and reiteration.
\newblock {\em Canad. J. Math.}, 52:920--960, 2000.

\bibitem[EOP02]{EOP02:RILF}
W.~D. Evans, B.~Opic, and L.~Pick.
\newblock Real {I}nterpolation with {L}ogarithmic {F}unctors.
\newblock {\em J. Inequal. Appl.}, 7(2):187--269, 2002.

\bibitem[FFG18]{FFG2018}
A.~Fiorenza, M.~R. Formica, and A.~Gogatishvili.
\newblock On grand and small {L}ebesgue and {S}obolev spaces and some
  applications to {PDE}'s.
\newblock {\em Differ. Equ. Appl.}, 10(1):21--46, 2018.

\bibitem[FK04]{FK2004}
A.~Fiorenza and G.~E. Karadzhov.
\newblock Grand and small {L}ebesgue spaces and their analogs.
\newblock {\em Z. Anal. Anwendungen}, 23(4):657--681, 2004.

\bibitem[FL06]{FaLe06}
W.~Farkas and H.-G. Leopold.
\newblock Characterisations of function spaces of generalised smoothness.
\newblock {\em Ann. Mat. Pura Appl. (4)}, 185(1):1--62, 2006.

\bibitem[GK03]{GoKe}
M.~L. Gol{\cprime}dman and R.~A. Kerman.
\newblock On the optimal embedding of {C}alder\'{o}n spaces and of generalized
  {B}esov spaces.
\newblock {\em Tr. Mat. Inst. Steklova}, 243:161--193, 2003.
\newblock (Russian) Translated in Proc. Steklov Inst. Math. {243}: 154--184,
  2003.

\bibitem[GNO10]{GNO10:PotAnal}
A.~Gogatishvili, J.~S. Neves, and B.~Opic.
\newblock Optimal embeddings of {B}essel-potential-type spaces into generalized
  {H}\"older spaces involving $k$-modulus of smoothness.
\newblock {\em Potential Anal.}, 32(3):201--228, 2010.

\bibitem[GO05]{GuOp04}
P.~Gurka and B.~Opic.
\newblock Sharp embeddings of {B}esov spaces with logarithmic smoothness.
\newblock {\em Rev. Mat. Complut.}, 18(1):81--110, 2005.

\bibitem[GO07]{GO07:SEBtS}
P.~Gurka and B.~Opic.
\newblock Sharp embeddings of {B}esov-type spaces.
\newblock {\em J. Comput. Appl. Math.}, 208:235--269, 2007.

\bibitem[Gol87]{Gol-Russian}
M.~L. Gol{\cprime}dman.
\newblock On the imbedding of a {N}ikol'skij--{B}esov space in a weighted
  {L}orentz space.
\newblock {\em Trudy Mat. Inst. Steklov.}, 180:93--95, 1987.
\newblock (Russian).

\bibitem[Gol07]{Gol07}
M.~L. Gol{\cprime}dman.
\newblock Rearrangement invariant envelopes of generalized {B}esov, {S}obolev,
  and {C}alder\'on spaces.
\newblock In T.~Iwaniec V.~I.~Burenkov and S.~K. Vodopyanov, editors, {\em The
  interaction of analysis and geometry}, volume 424 of {\em Contemp. Math.},
  pages 53--81. Amer. Math. Soc., Providence, RI, 2007.

\bibitem[GOT05]{GOT2002Lrrinwsvf}
A.~Gogatishvili, B.~Opic, and W.~Trebels.
\newblock Limiting reiteration for real interpolation with slowly varying
  functions.
\newblock {\em Math. Nachr.}, 278:86--107, 2005.

\bibitem[GOTT14]{GOTT2014}
A.~Gogatishvili, B.~Opic, S.~Tikhonov, and W.~Trebels.
\newblock Ulyanov-type inequalities between {L}orentz-{Z}ygmund spaces.
\newblock {\em J. Fourier Anal. Appl.}, 20(5):1020--1049, 2014.

\bibitem[HS93]{H-Step}
H.~P. Heinig and V.~D. Stepanov.
\newblock Weighted {H}ardy inequalities for increasing functions.
\newblock {\em Can. J. Math.}, 45:104--116, 1993.

\bibitem[KL87]{KL87}
G.~A. Kaljabin and P.~I. Lizorkin.
\newblock Spaces of functions of generalized smoothness.
\newblock {\em Math. Nachr.}, 133:7--32, 1987.

\bibitem[Lai93]{Lai:1993wt}
S.~Lai.
\newblock {Weighted Norm Inequalities for General Operators on Monotone
  Functions}.
\newblock {\em Trans. Amer. Math. Soc.}, 340(2):811--836, 1993.

\bibitem[Mar00]{VojislavMaric:RVDE00}
V.~Mari\'c.
\newblock {\em Regular Variation and Differential Equations}, volume 1726 of
  {\em Lecture Notes in Mathematics}.
\newblock Springer-Verlag, Berlin, 2000.

\bibitem[Mou01]{Moura}
S.~D. Moura.
\newblock Function spaces of generalised smoothness.
\newblock {\em Dissertationes Math. (Rozprawy Mat.)}, 398:88, 2001.

\bibitem[Net87]{Net89}
Yu.~V. Netrusov.
\newblock Theorems for the embedding of {B}esov spaces into ideal spaces.
\newblock {\em Zap. Nauchn. Sem. Leningrad. Otdel. Mat. Inst. Steklov. (LOMI)},
  159:69--82, 1987.
\newblock (Russian) Translated in J. Soviet Math. 47 (6): 2871--2881, 1989.

\bibitem[Nev02]{Nev02:LKSBRPE}
J.~S. Neves.
\newblock Lorentz-{K}aramata spaces, {B}essel and {R}iesz potentials and
  embeddings.
\newblock {\em Dissertationes Math. (Rozprawy Mat.)}, 405:46 pp., 2002.

\bibitem[OK90]{OK90:HTI}
B.~Opic and A.~Kufner.
\newblock {\em Hardy-type inequalities}.
\newblock Pitman Research Notes in Math. Series 219, Longman Sci. \& Tech.,
  Harlow, 1990.

\bibitem[Saw90]{Swa90:BCOCLS}
E.~T. Sawyer.
\newblock Boundedness of classical operators on classical {L}orentz spaces.
\newblock {\em Studia Math.}, 96:145--158, 1990.

\bibitem[ST95]{SiTr}
W.~Sickel and H.~Triebel.
\newblock H\"{o}lder inequalities and sharp embeddings in function spaces of
  {$B^s_{pq}$} and {$F^s_{pq}$} type.
\newblock {\em Z. Anal. Anwendungen}, 14(1):105--140, 1995.

\bibitem[Tri06]{Tri06}
H.~Triebel.
\newblock {\em Theory of function spaces. {III}}, volume 100 of {\em Monographs
  in Mathematics}.
\newblock Birkh\H{a}user Verlag, Basel, 2006.

\bibitem[Zyg57]{Zyg57:TS}
A.~Zygmund.
\newblock {\em Trigonometric Series}, volume~I.
\newblock Cambridge University Press, Cambridge, 1957.

\end{thebibliography}
\def\cprime{$'$}

\end{document}